# Stochastic optimization of a mixed moving average process for controlling non-Markovian streamflow environments

Hidekazu Yoshioka[1, *], Tomohiro Tanaka[2], Yumi Yoshioka[1], Ayumi Hashiguchi[1]

[1] Shimane University, Nishikawatsu-cho 1060, Matsue, Shimane 690-8504, Japan

[2] Kyoto University, Kyotodaigaku-katsura, Nishikyo-ku, Kyoto, Kyoto 615-8540, Japan

* Corresponding author: E-mail: yoshih@life.shimane-u.ac.jp

**Abstract**

We investigated a cost-constrained static ergodic control problem of the variance of measure-valued affine processes and its application in streamflow management. The controlled system is a jump-driven mixed moving average process that generates realistic subexponential autocorrelation functions, and the "static" nature of the control originates from a realistic observability assumption in the system. The Markovian lift was effectively used to discretize the system into a finite-dimensional process, which is easier to analyze. The resolution of the problem is based on backward Kolmogorov equations and a quadratic solution ansatz. The control problem has a closed-form solution, and the variance has both strict upper and lower bounds, indicating that the variance cannot take an arbitrary value even when it is subject to a high control cost. The correspondence between the discretized system based on the Markovian lift and the original infinite-dimensional one is discussed. Then, a convergent Markovian lift is presented to approximate the infinite-dimensional system. Finally, the control problem was applied to real cases using available data for a river reach. An extended problem subject to an additional constraint on maintaining the flow variability was also analyzed without significantly degrading the tractability of the proposed framework.

## 1. Introduction

### 1.1 Background

Streamflow as a stochastic process provides dynamic habitats for several aquatic fauna [1, 2] and determines the transport phenomena of sediment and water-quality indices [3, 4]. Understanding its stochastic nature as a hydrological driver is also important for human lives, including in disaster risk management [5, 6] and water supply [7, 8].

Some mathematical models for the stochastic streamflow, including the water level and discharge at a fixed observation station, are jump-driven stochastic differential equations (SDEs) containing flood and recession terms [9]. The recession mechanism is assumed to be an exponential decay in continuous-time models [10, 11], which appears as single-step linear decay in discrete-time models [12]. These models have been employed in both theoretical analysis and engineering applications owing to their high compatibility; thus, the statistics of key hydrological variables for evaluating the stochastic streamflow are available in closed forms. Such examples include but are not limited to the active streamflow length [13], fluvial erosion [14], and water quality assessment [15].

Conventional jump-driven SDEs and their discrete analogs employ exponential autocorrelation functions (ACFs), whereas real hydrological data do not necessarily satisfy this property. Subexponential ACFs have been reported for discharge time series in mountainous river environments [16, 17], highly regulated rivers [18, 19], and the turbidity time series in rivers flowing near erodible areas [20]. Reproducing the memory of streamflow accurately is important in theory to better understand their dynamics [21] as well as in applications because the memory properties critically affect hydrological forecasting [22].

In practice, models that capture subexponential memory should be tractable. (Conditional) expectations of some objective functions should be computable analytically or efficiently by numerical computation as many engineering problems involving river environments are stochastic optimization problems. Non-Markov stochastic process models, including Volterra processes containing the memory terms [23] and mixed moving average (MMA) processes as multiple timescale models [24, 25], are candidates in economics, where the stochastic optimization of advanced stochastic processes have long been extensively studied. Recently, these models have been employed in discharge time series [17, 26], where a superposition of Ornstein–Uhlenbeck processes (a supOU process) as an affine MMA process was utilized. The largest advantage of the MMA process would be its ability to model subexponential ACFs that cannot be reproduced by common SDE models like Ornstein–Uhlenbeck processes having exponential ACFs. The first author reported that the polynomial decay of the discharge times series in mountainous river environments could be accurately captured through supOU processes [17]. Although the above-mentioned models serve as useful stochastic processes having a long memory, there are still mathematical and/or computational difficulties due to the inherent non-Markovian properties. These difficulties have partly been resolved using the Markovian lift to reconsider a non-Markov process as a system of measure-valued Markov processes [27], the latter being more amenable to analyses. However, the role and effectiveness of the Markovian lift in modelling streamflow dynamics have not been well-studied.

Another issue to be considered in the application of non-Markov processes, especially in their optimal control, is the controller design. In most basic linear–quadratic (LQ) control problems, the controller fully depends on the lifted Markovian variables for MMA and stochastic Volterra processes [26, 28], whereas it is difficult to implement it in real applications. Indeed, river discharge governed by an MMA process is an observable variable, whereas all measure-valued processes as elements constituting the discharge cannot be observed even though they are mathematically well-defined. Namely, an integral of the measure-valued processes can be observed but not each of them. Therefore, the controller should be designed based on only observable variables, resulting in a nonstandard stochastic optimization problem with a static feedback controller [29]. Typically, constraining the controller design prevents the application of the standard approach in solving the Riccati equation and further leads to NP-hardness [30] and the existence of many local optimizers [31]; however, these issues can be resolved if the problem admits a closed-form solution verified to be optimal. This is the motivation behind this study.

### 1.2 Objective and contribution

In this study, we investigated the application of an MMA process to a static feedback control problem of river discharge as water resources. The system to be controlled is a jump-driven affine SDE based on an MMA process. The controller contains several real parameters with which the average discharge can be designed to be an arbitrary value, and the variance of the discharge is controllable subject to the control cost. We consider a minimization problem of variance of the difference between the discharge and its prescribed target value subject to an expectation constraint of a quadratic control cost. The problem is fairly simple but has not been reported in the literature as it is not a classical LQ control problem with full observation. Our assumption on the controller agrees well with the fact that the discharge given as integration of measure-valued processes is observable, whereas each of the latter processes is not.

The system dynamics are infinite-dimensional since the MMA process is an integral of independent measure-valued processes, posing a theoretical challenge. Thus, the Markovian lift [32] is applied to the MMA process to derive a finite-dimensional counterpart. With this finite-dimensional system, the conventional theory of SDEs [33] could be applied to our problem, and the Markovian lift is justified because the finite-dimensional MMA process converges to the original infinite-dimensional one in the sense of distributions [34].

In this study, first, we discretize the system using a Markovian lift, solve the finite-dimensional optimization problem, and then take its infinite-dimensional limit (**Figure 1**). Using proper backward Kolmogorov equations (BKEs), we show that both the controlled variance and cost are found analytically once a controller is specified. This remarkable tractability of the proposed model allows for the efficient analysis of the optimal controller and the optimized variance and cost. As a byproduct, we obtain that the variance cannot take an arbitrary value even subject to a huge control cost. We also discuss the solution procedure for the optimization problem without bypassing the Markovian lift. This procedure employs a nonlocal partial differential equation (PDE) in an infinite dimension and is, therefore, a difficult mathematical subject, whereas the solution can be guessed as in the finite-dimensional case. The approach based on BKE can also be useful for tackling more complicated problems in the future. This study serves as a basis for such cases. Furthermore, we present a provably-convergent Markovian lift with quantile-based discretization. Both mathematical analysis and optimal control or optimization of MMA processes have not been extensively studied so far in hydrology and related fields.

Finally, the proposed model was employed in the modelling and hypothetical optimization of river discharge in an existing river environment in Japan. The target site is a river reach where the coexistence between the environment and humanities is of current concern. The stochastic flow field in the river reach driven by the upstream boundary conditions and lateral inflows were emulated using a verified hydrodynamic model. A similar and slightly simpler approach for generating a random field has been employed in the context of electricity supply [35], whereas the stochastic processes and control objectives are different. The MMA process for capturing polynomially-decaying ACF was identified at each point in the river reach using the generated flow field, and the sensitivity of the control problem was analyzed. An extended problem subject to an additional constraint on maintaining the flow variability was also analyzed without critically degrading the tractability of the proposed mathematical framework. This study contributes to the new modelling and analysis of river discharge based on the MMA process and static control.

The rest of the paper is structured as follows. **Section 2** explains our mathematical model based on the supCBI process and presents the Markovian lift. The controller design is also explained in this section. **Section 3** introduces our optimization problem and presents several mathematical results, including the main results **Propositions 7 and 8** for finite- and infinite-dimensional cases (**Sections 3.1** to **3.4**) and a convergence argument of the Markovian lift (**Section 3.5**) that supports our computational results. A solution algorithm for solving the problem, which is based on the mathematical analysis results, is also explained in this section **(Section 3.6**). **Section 4** studies its application to an existing river environment. **Section 5** gives a summary and future perspectives of this study. **Appendices** in the supplementary material contains proofs of propositions and supporting data concerning the computation and application.

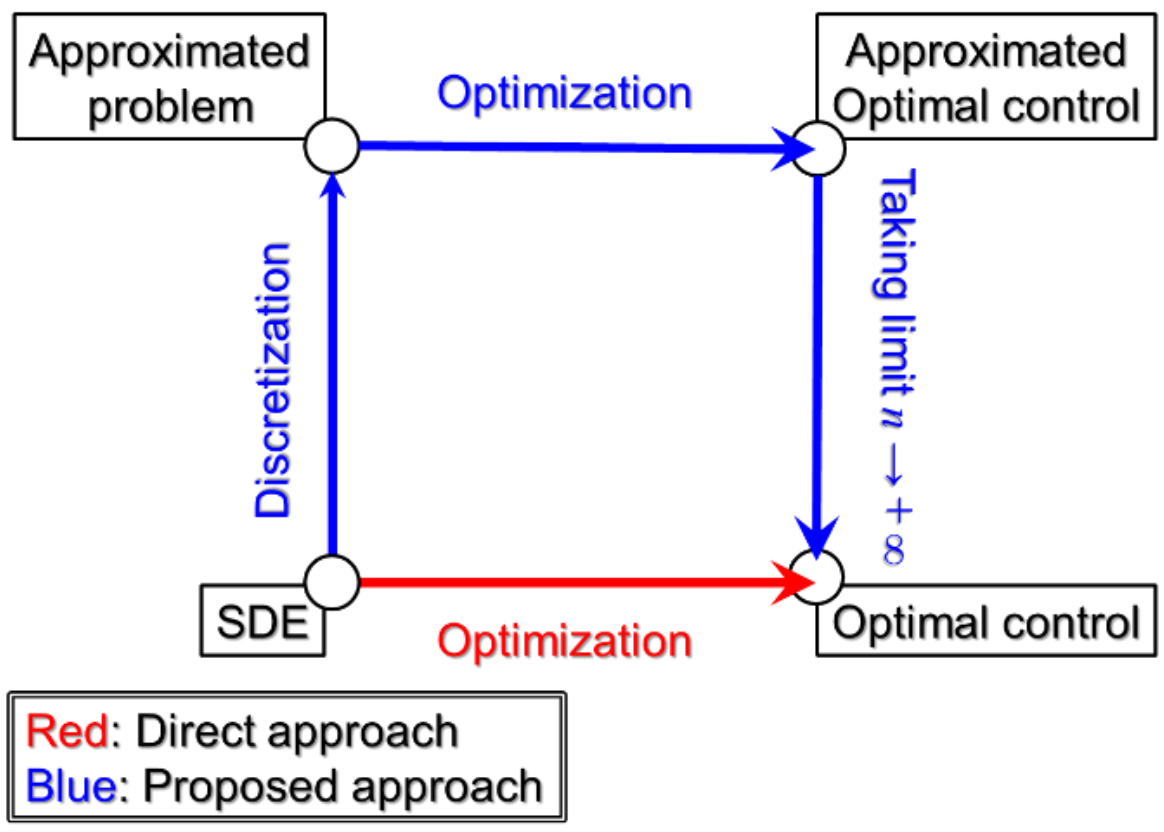

**Figure 1.** A schematic diagram of the proposed and direct approaches.

## 2. Mathematical model

### 2.1 Discharge

We consider an infinite time horizon $\mathbb{R}$ for a time-average, ergodic control problem. The discharge at a fixed location in a longitudinal river is considered a continuous-time scalar variable $Q$, which is denoted as $Q_t$ at time $t$. The discharge is assumed to be the sum of the inflows from the upstream $\underline{Q}+Y_t$ with a constant baseflow $\underline{Q} \geq 0$, the jump-driven runoff process $Y=(Y_t)_{\in\mathbb{R}}$ as a given uncontrollable process, and the water addition/abstraction $C=(C_t)_{\in\mathbb{R}}$ as a function of the observable $Q$:

$$Q_t = C_t + \underline{Q} + Y_t, \quad t \in \mathbb{R}. \tag{1}$$

This equation expresses the mass conservation law to be satisfied at the point of interest in the river.

The right-hand side of (1) is specified as follows. Our formulation is based on measure-valued processes, whereas the optimization result is simple and intuitive (**Section 3**). We set $Y$ as a superposition of continuous-state branching processes with immigration (a supCBI process), which is an MMA process to describe the river discharge time series having a subexponential ACF [32]:

$$Y_t = \int_0^{+\infty} Y_t^{(r)}(\mathrm{d}r), \quad t > 0 \tag{2}$$

where $Y_t^{(r)}(\mathrm{d}r)$ is a measure-valued process given for $t \in \mathbb{R}$ and $r > 0$, which is given by

$$Y_t^{(r)}(\mathrm{d}r) = \int_{-\infty}^{t} e^{-r(t-s)} \int_0^{A\pi(\mathrm{d}r)+rBY_{t-}^{(r)}(\mathrm{d}r)} \int_0^{\infty} z\mu^{(r)}(\mathrm{d}w,\mathrm{d}z,\mathrm{d}s) \tag{3}$$

or in a differential form

$$\mathrm{d}Y_t^{(r)}(\mathrm{d}r) = -rY_t^{(r)}(\mathrm{d}r)\mathrm{d}t + \int_0^{A\pi(\mathrm{d}r)+rBY_{t-}^{(r)}(\mathrm{d}r)} \int_0^{\infty} z\mu^{(r)}(\mathrm{d}w,\mathrm{d}z,\mathrm{d}t), \tag{4}$$

with parameters $A>0, B\geq 0$, a probability measure $\pi(\mathrm{d}r)$ of positive random variables satisfying $\int_0^{+\infty} r^{-1}\pi(\mathrm{d}r) < +\infty$ and a Poisson random measure $\mu^{(r)}$ on $\mathbb{R}^3$ having the intensity measure $\mathrm{d}w \times \nu(\mathrm{d}z) \times \mathrm{d}t$ with a positive-jump background driving Lévy measure $\nu(\mathrm{d}z)$ having first and second moments; that is,

$$\int_{(-\infty,0]} \nu(\mathrm{d}z) = 0, \quad \int_{(0,+\infty)} \min\{1,z\}\,\nu(\mathrm{d}z) < +\infty, \text{ and } M_k = \int_{(0,+\infty)} z^k \nu(\mathrm{d}z) < +\infty \quad (k=1,2). \tag{5}$$

The Poisson random measure $\mu^{(r)}$ is defined in the real time real line $\mathbb{R}$ based on [36], so that we can consider a stationary stochastic process at any finite time $t$. Here, notice that in (3) the integration range of the time starts from $-\infty$. The first and second conditions of (5) ensure that the background driving Lévy measure has a finite variation with positive jumps, whereas the last one is a technical assumption that well-poses our optimization problem. This condition is not restrictive in practice (**Section 4**).

The coefficient $A\pi(\mathrm{d}r) + rBY_{t-}^{(r)}(\mathrm{d}r)$ is a proportional coefficient of the jump rate of $Y_t^{(r)}(\mathrm{d}r)$, where the first and second terms represent the constant and self-excited parts corresponding to clustered flood events, respectively. We assume $B$ is small such that $BM_1 < 1$, indicating that self-excited jumps (clustered flood events) are not extremely frequent in the discharge. This condition is necessary and sufficient for the stationarity and existence of the supCBI process [32]. If $B = 0$, then the self-excited part is dropped out, and the supCBI process reduces to the superposition of the Ornstein–Uhlenbeck processes.

Equations (3) and (4) show that inflow (2) is a superposition of measure-valued processes having different values of the reversion rate $r > 0$, which represents a situation where river discharge is considered a scalar variable containing both fast (surface) and slow (underground) runoff processes having multiple time scales [37, 38]. The coexistence of multiple time scales is represented by the probability measure $\pi$, which can yield the subexponential ACF of $Y$ if specified adequately. The stationary ACF of $Y$ with a time lag $\tau \geq 0$ is given as follows [32]:

$$\mathrm{ACF}(\tau) = \left( \int_0^{+\infty} r^{-1} \pi(\mathrm{d}r) \right)^{-1} \int_0^{+\infty} r^{-1} e^{-D\tau r} \pi(\mathrm{d}r) \tag{6}$$

where $D = 1 - BM_1 \in (0,1)$. For example, choosing a Gamma one $\pi(\mathrm{d}r) \sim r^{\alpha-1} e^{-\frac{r}{\beta}} \mathrm{d}r$ ($\alpha > 1, \beta > 0$) yields $\mathrm{ACF}(\tau) = (1 + \beta\tau)^{-(\alpha-1)}$, which decays polynomially.

## 2.2 Markovian lift

Herein, we approximate the supCBI process using a Markovian lift. We set the degree of freedom $n \in \mathbb{N}$, a positive and strictly increasing sequence $\{r_i\}_{1 \leq i \leq n}$ of the discretized reversion rate, and a positive sequence $\{c_i\}_{1 \leq i \leq n}$ with $\sum_{i=1}^{n} c_i = 1$ as discrete probabilities. Then, the *n*-dimensional supCBI process $\left\{ Y^{(r_i)} \right\}_{1 \leq i \leq n}$ is set as a unique càdlàg solution to the system

$$\mathrm{d}Y_t^{(r_i)} = -\rho_i Y_t^{(r_i)} \mathrm{d}t + \int_0^{c_i A + \rho_i B Y_{t-}^{(r_i)}} \int_0^{\infty} z \mu_i(\mathrm{d}u, \mathrm{d}z, \mathrm{d}t), \quad t \in \mathbb{R}, \quad 1 \leq i \leq n. \tag{7}$$

Each $\mu_i$ ($1 \leq i \leq n$) is a mutually independent Poisson random measure with a compensator $\mathrm{d}u \times v(\mathrm{d}z) \times \mathrm{d}t$. Then, the finite-dimensional counterpart of (1) is set as

$$Q_t = C_t + \underline{Q} + \sum_{i=1}^{n} Y_t^{(r_i)} \equiv C_t + \underline{Q} + Y_{n,t}, \quad t \in \mathbb{R}. \tag{8}$$

The Markovian lift, which is essentially a space discretization of the supCBI process in phase space, serves well such that $Y_n$ converges to $Y$ in the sense of distribution if the parameters $\{r_i\}_{1 \leq i \leq n}$ and $\{c_i\}_{1 \leq i \leq n}$ are chosen adequately. We present an example in **Section 3**. For the meantime, we assume such a discretization in a Markovian lift.

## 2.3 Controller design

The control $C$ is specified as a function of the discharge $Q$ observable continuously in time. We consider the following integral form:

$$C_t = u\int_{-\infty}^{t} e^{-\rho(t-s)}\left(Q_s + W\right)\mathrm{d}s\ ,\quad t\in\mathbb{R}\ . \tag{9}$$

This controller is determined by $u\in\mathbb{R}$, where $u>0$ (resp., $u<0$) is the positive (resp., negative) feedback, the memory effect is modulated by $\rho>0$, where a larger $\rho$ is a more myopic controller, and the constant shift $W\in\mathbb{R}$. We assume $\rho>u$. The system, even without noise, becomes exponentially unstable if $\rho\le u$.

Substituting (9) into (8), we obtain the following equation:

$$Q_t = u\int_{-\infty}^{t} e^{-\rho(t-s)}\left(Q_s + W\right)\mathrm{d}s + \underline{Q} + \sum_{i=1}^{n} Y_t^{(r_i)}\ ,\quad t\in\mathbb{R}\ . \tag{10}$$

The shifted process $X=\left(X_t\right)_{t\in\mathbb{R}}$ is defined as

$$X_t = Q_t - \frac{\rho}{\rho-u}\left(\underline{Q} + \frac{uW}{\rho}\right),\quad t\in\mathbb{R}\ . \tag{11}$$

Combining (10) and (11) yields

$$X_t = u\int_{-\infty}^{t} e^{-\rho(t-s)} X_s \mathrm{d}s + \sum_{i=1}^{n} Y_t^{(r_i)}\ ,\quad t\in\mathbb{R}\ . \tag{12}$$

The shifted process $X$ is still observable given $\underline{Q},\rho,u,W$.

From (12), we obtain the SDE

$$\mathrm{d}X_t = \left(uX_t - u\rho\int_{-\infty}^{t} e^{-\rho(t-s)} X_s \mathrm{d}s\right)\mathrm{d}t + \sum_{i=1}^{n}\mathrm{d}Y_t^{(r_i)} = c_t\mathrm{d}t + \sum_{i=1}^{n}\mathrm{d}Y_t^{(r_i)}\ ,\quad t\in\mathbb{R} \tag{13}$$

with the control rate

$$c_t = -\left(\rho-u\right)X_t + \rho\sum_{i=1}^{n} Y_t^{(r_i)}\ . \tag{14}$$

The SDE representation (13) shows that the first term on the right-hand side of (13) is the control effect, and the second term is an increment of the inflow as exogenic noise.

### 2.4 Optimization problem

The optimization problem here is as follows:

$$\text{find } \inf_{\rho,u,W} J\ \text{, where } J = \limsup_{T\to+\infty}\frac{1}{T}\mathbb{E}\left[\int_0^T\left(Q_t-\hat{Q}\right)^2\mathrm{d}t\right] \tag{15}$$

subject to the cost constraint

$$K = \limsup_{T\to+\infty}\frac{1}{T}\mathbb{E}\left[\int_0^T c_t^2\mathrm{d}t\right]\le\bar{K} \tag{16}$$

and the parameter constraints

$$\rho,u,W\in\mathbb{R}\ ,\quad \rho>0\ ,\quad \rho>u\ , \tag{17}$$

and

$$\mathbb{E}\left[X_t\right]=\frac{\rho}{\rho-u}\mathbb{E}\left[Y_{n,t}\right]. \tag{18}$$

Here, $\hat{Q}>0$ is a prescribed target discharge, and $\bar{K}>0$ is a prescribed upper bound of the control cost evaluated as a time-average integral. The objective expressed by (15) implies that the deviation of the discharge from the target is penalized quadratically, and this deviation is evaluated by averaging the time. Therefore, this optimization problem is a long-run optimal control problem. The constraint (16) means that the controlling effort of the streamflow dynamics is strictly bounded from above. The constraint in (18) is to achieve the target discharge in the mean as $\mathbb{E}\left[Q_t\right]=\hat{Q}$ ($\mathbb{E}\left[X_t\right]=\hat{X}$). This is a natural constraint because the discharge should be controlled at or around the target $\hat{Q}$. By rewriting $\mathbb{E}\left[Q_t\right]=\hat{Q}$, we obtain (18), where $\mathbb{E}\left[Y_{n,t}\right]=\frac{AM_1}{D}\sum_{i=1}^{n}\frac{c_i}{r_i}$ [32].

Using the shifted process $X$ in (11), we can rewrite (15) as

$$J=\limsup_{T\to+\infty}\frac{1}{T}\mathbb{E}\left[\int_0^T\left(X_t-\hat{X}\right)^2\mathrm{d}t\right] \text{ with } \hat{X}=\hat{Q}-\frac{\rho}{\rho-u}\left(\underline{Q}+\frac{uW}{\rho}\right). \tag{19}$$

Our optimization problem is seen an LQ control problem but is not a classical one whose resolution reduces to finding an appropriate solution to a Riccati equation [39]. This is because not all the state variables are observable but the discharge $Q$ (or $X$) is, implying that we cannot resort to dynamic programming as a well-stylized methodology for solving stochastic control problems. We resolve this issue by finding a closed-form solution to the optimization problem with a direct calculation, as shown in the next section.

## 3. Solution to the optimization problem

In this section, we show that the expectations in (15) and (16) can be explicitly computed. To this end, we derive, solve, and verify the BKEs of these integrals.

### 3.1 Backward Kolmogorov equation of the variance

Based on the SDE representations in (13) combined with (15), a classical Dynkin's formula gives [e.g., 15]

$$\begin{aligned}&-J+\left(x-\hat{X}\right)^2+\left\{-(\rho-u)x+\sum_{i=1}^{n}(\rho-r_i)y_i\right\}\frac{\partial\Phi}{\partial x}-\sum_{i=1}^{n}r_iy_i\frac{\partial\Phi}{\partial y_i}\\&+\sum_{i=1}^{n}\left(c_iA+r_iBy_i\right)\int_0^{+\infty}\left(\Phi\left(x+z,\ldots,y_i+z,\ldots\right)-\Phi\right)\nu\left(\mathrm{d}z\right)=0\end{aligned} \tag{20}$$

and the solution is a pair of the constant $J\in\mathbb{R}$ (this turns to be the objective (15)) and a smooth function $\Phi:\mathbb{R}^{n+1}\to\mathbb{R}$ growing at most polynomially. Here, $J$ is important in practice, and the function $\Phi$ is an auxiliary variable for computing $J$.

Our first main result is **Proposition 1**, which states that $J$ is found explicitly. This proposition is proved in **Appendix A.1** because it is technical, and $\Phi$ is also provided there. We set $q=\dfrac{\rho}{\rho-u}>0$. The key of the proof is a guessed-solution technique with a quadratic ansatz of $\Phi$.

***Proposition 1:***

$$
\begin{aligned}
J &= \left( \hat{X} - \frac{qAM_1}{D}\sum_{i=1}^{n}\frac{c_i}{r_i} \right)^2 + \frac{1}{2}\frac{AM_2}{D^2}\sum_{i=1}^{n}\frac{c_i}{r_i}\frac{r_iD+q^2(\rho-u)}{r_iD+\rho-u} \\
&= \left( \hat{X} - q\mathbb{E}[Y_n] \right)^2 + \mathrm{Var}[Y_n]\left( \sum_{i=1}^{n}\frac{c_i}{r_i} \right)^{-1}\sum_{i=1}^{n}\frac{c_i}{r_i}\frac{r_iD+q^2(\rho-u)}{r_iD+\rho-u}
\end{aligned}
. \tag{21}
$$

**Proposition 1** suggests that the objective consists of two terms: the first term represents the deviation between the target $\hat{X}$ and the scaled mean inflow $\mathbb{E}[Y_n]=\dfrac{AM_1}{D}\sum_{i=1}^{n}\dfrac{c_i}{r_i}$, and the second term is proportional to the variance $\mathrm{Var}[Y_n]=\dfrac{1}{2}\dfrac{AM_2}{D^2}\sum_{i=1}^{n}\dfrac{c_i}{r_i}$ of the inflow.

### 3.2 Backward Kolmogorov equation of the cost

As in the previous subsection, we analytically evaluate the cost $K$ by solving BKE. First, we have

$$
\begin{aligned}
K &= \limsup_{T\to+\infty}\frac{1}{T}\mathbb{E}\left[ \int_0^T c_t^2 \mathrm{d}t \right] \\
&= \limsup_{T\to+\infty}\frac{1}{T}\mathbb{E}\left[ \int_0^T \left( -(\rho-u)X_t + \rho\sum_{i=1}^{n}Y_t^{(r_i)} \right)^2 \mathrm{d}t \right] \\
&= (\rho-u)^2\limsup_{T\to+\infty}\frac{1}{T}\mathbb{E}\left[ \int_0^T \left( X_t - \frac{\rho}{\rho-u}\sum_{i=1}^{n}Y_t^{(r_i)} \right)^2 \mathrm{d}t \right] \\
&\equiv (\rho-u)^2 L
\end{aligned}
. \tag{22}
$$

Therefore, we should find $L$. Again using Dynkin's formula, we obtain BKE

$$
\begin{aligned}
&-L + \left( x - q\sum_{i=1}^{n} y_i \right)^2 + \left\{ -(\rho-u)x + \sum_{i=1}^{n}(\rho-r_i)y_i \right\}\frac{\partial\Psi}{\partial x} - \sum_{i=1}^{n} r_i y_i \frac{\partial\Psi}{\partial y_i} \\
&+ \sum_{i=1}^{n}(c_iA + r_iBy_i)\int_0^{+\infty}\left( \Psi(x+z,\ldots,y_i+z,\ldots) - \Psi \right)v(\mathrm{d}z) = 0
\end{aligned}
\tag{23}
$$

whose solution is a pair of the constant $L\in\mathbb{R}$ and a smooth function $\Psi:\mathbb{R}^{n+1}\to\mathbb{R}$ growing at most polynomially. Also, $\Psi$ serves as an auxiliary variable for computing $J$.

Our second main result is **Proposition 2**, which states that $L$ and, hence, $K$ is obtained explicitly. This proposition is proved in **Appendix A.2** because it is technical, and $\Psi$ is also provided. The proof is also based on a guessed solution.

***Proposition 2:***

$$L=\frac{1}{2}(1-q)^2\frac{AM_2}{D^2}\sum_{i=1}^{n}\frac{c_i}{r_i}\frac{r_iD}{r_iD+\rho-u} \tag{24}$$

*and hence*

$$\begin{aligned}K&=\frac{1}{2}(\rho-u)^2(1-q)^2\frac{AM_2}{D^2}\sum_{i=1}^{n}\frac{c_i}{r_i}\frac{r_iD}{r_iD+\rho-u}\\&=(\rho-u)^2(1-q)^2\operatorname{Var}[Y_n]\left(\sum_{i=1}^{n}\frac{c_i}{r_i}\right)^{-1}\sum_{i=1}^{n}\frac{c_i}{r_i}\frac{r_iD}{r_iD+\rho-u}\end{aligned}. \tag{25}$$

**Proposition 2** is further explained in the next subsection.

### 3.3 Closed-form solution

**Propositions 2 and 3** allow us to obtain a closed-form solution to the minimization problem of $J$. From (18), we obtain

$$q=\frac{D\hat{X}}{AM_1}\sum_{i=1}^{n}\frac{c_i}{r_i}=\frac{\mathbb{E}[X_t]}{\mathbb{E}[Y_{n,t}]}=\frac{\hat{X}}{\mathbb{E}[Y_{n,t}]}. \tag{26}$$

The right-hand side of (26) is a constant as it is free from $\rho,u,W$. Now, we can reduce **Proposition 1** as the first term of (21) disappears in (26).

***Proposition 3:*** *Under (26), it follows that*

$$J=\operatorname{Var}[Y_n]\left(\sum_{i=1}^{n}\frac{c_i}{r_i}\right)^{-1}\sum_{i=1}^{n}\frac{c_i}{r_i}\frac{r_iD+q^2(\rho-u)}{r_iD+\rho-u}. \tag{27}$$

Furthermore, the constraint (26) explains the coefficient $(1-q)^2$ appearing in $L$.

***Proposition 4:*** *Under (26) with $q=1$, the target discharge equals the average of the inflow, and hence $K=0$.*

From (26), we have $u=0$ and, hence, $C_t=0$ all the time. Given **Proposition 4**, (25) is reasonable as there are some control costs if the target discharge equals the average of the inflow. Then, we solve the optimization problem by invoking (26).

***Remark 1:*** By **Propositions 2 and 3**: the parameter $W$ does not play an essential role in (26). We therefore set $W=0$ as it is a redundant parameter. Notably, this redundancy was initially nontrivial.

We investigated $J$ without the cost constraint (16). **Proposition 5** can be checked directly.

***Proposition 5:*** *Fix* $q>0$. *Consider* $J$ *as a function of* $h\equiv\rho-u\ge 0$. *Then, it follows that*

$$\mathrm{Var}\left[Y_n\right]\min\left\{1,q^2\right\}\le J\left(h\right)\le \mathrm{Var}\left[Y_n\right]\max\left\{1,q^2\right\},\quad h\ge 0. \tag{28}$$

*Therefore,* $J$ *is essentially bounded both from below and above. In addition,* $J\left(h\right)$ *strictly increases (resp., decreasing) for* $h\ge 0$ *if* $q>1$ *(resp.,* $0<q<1$*) and equals* $\mathrm{Var}\left[Y_n\right]$ *if* $q=1$.

**Proposition 5** analytically characterizes the dependence of the cost on the parameters $q$ and $\rho-u$. It is essentially bounded once $q$ is given, implying that the variance cannot be arbitrarily decreased or increased by designing the controller.

***Remark 2:*** Our problem is a minimization problem of $J$. However, some applied problems may optimize $J$ at or around a prescribed level as the natural flow regime is typically not constant in time, and there would be an optimal disturbance level to maintain the food web structure and biodiversity [40]. **Proposition 5** is promising in this view because such a problem becomes infeasible with our controller. The problem may be overcome using a more complex and, possibly nonlinear controller instead of (9), with which the problem resolution requires a more complicated approach .

Concerning the cost $K$, we obtain the following proposition that can be checked directly.

***Proposition 6:*** *Fix* $q>0$ *with* $q\neq 1$. *Consider* $K$ *as a function of* $h\equiv\rho-u\ge 0$. *Then, we have that* $K\left(0\right)=0$, $\lim_{h\to+\infty}K\left(h\right)=+\infty$, *and* $K\left(h\right)$ *strictly increases for* $h\ge 0$. *Therefore, given* $\bar{K}>0$, *there exists a unique positive solution* $h=\bar{h}$ *such that*

$$K\left(h\right)=\left(1-q\right)^2\mathrm{Var}\left[Y_n\right]\left(\sum_{i=1}^{n}\frac{c_i}{r_i}\right)^{-1}\bar{h}^2\sum_{i=1}^{n}\frac{c_i}{r_i}\frac{r_iD}{r_iD+\bar{h}}=\bar{K}. \tag{29}$$

In practice, the solution $h=\bar{h}$ to (29) can be obtained by using a common numerical method. We will use the Picard method.

Based on **Propositions 5 and 6**, we can solve the optimization problem in a closed form.

***Proposition 7:*** *The minimization problem of* $J$ *subject to (16) and (17) under (26) is solved as follows:*

*(a) If* $q=1$, *then we have* $\min J=\mathrm{Var}\left[Y_n\right]$. *The minimizing couple* $\left(\rho,u\right)$ *is* $u=0$ *and arbitrary* $\rho>0$.

*(b) If* $q>1$ *(**Water Adding case**), then we have* $\inf_{h=\rho-u}J\left(h\right)=J\left(0\right)=\mathrm{Var}\left[Y_n\right]$. *The minimizing couple does not exist, but we have* $\lim_{h\to+0}J\left(h\right)=\mathrm{Var}\left[Y_n\right]$.

*(c) If* $0<q<1$ ***(Water Abstracting case)****, then we have* $\inf_{h=\rho-u} J(h)=J(\bar{h})=\mathrm{Var}[Y_n]$ *with* $\bar{h}$ *uniquely determined by (29). The minimizing couple* $(\rho,u)$ *satisfies* $q=\dfrac{\rho}{\rho-u}$ *and* $\rho-u=\bar{h}$*, and the minimum of* $J$ *is achieved by* $(\rho,u)=(q\bar{h},-(1-q)\bar{h})$.

The cost constraint in (16) works effectively only in the Water Abstracting case, in which the optimization problem is solved using the feasible pair $(\rho,u)=(q\bar{h},-(1-q)\bar{h})$. In the Water Adding case, no controller achieves infimum, although there are infinitely many controllers that can achieve the objective arbitrary closer to being optimal.

***Remark 3:*** In the formulation explained so far, the target is the discharge but not the water to be abstracted or added. However, the latter can also be used in the optimization problem. Let $Q_{\text{abs}}$ be the target abstracted discharge, such as the demand for irrigation or hydropower generation. If $q\in(0,1)$, then $\hat{Q}+Q_{\text{abs}}=\underline{Q}+\mathbb{E}[Y_n]$, or $\hat{Q}=\underline{Q}+\mathbb{E}[Y_n]-Q_{\text{abs}}$ at a stationary state. Considering $Q_{\text{abs}}$ instead of $\hat{Q}$, then this relationship can be used. The same applies to the infinite-dimensional case discussed below. The Water Adding case, although seems to be less realistic, can also be reexpressed in terms of the amount of water to be added.

### 3.4 Infinite-dimensional case

The minimization problem of the infinite-dimensional case can be formulated exactly as in (15)–(17) and the system dynamics

$$X_t=u\int_{-\infty}^{t}e^{-\rho(t-s)}X_s+Y_t,\quad t\in\mathbb{R} \tag{30}$$

or

$$\mathrm{d}X_t=\left(uX_t-u\rho\int_{-\infty}^{t}e^{-\rho(t-s)}X_s\right)\mathrm{d}t+\mathrm{d}Y_t=c_t\mathrm{d}t+\mathrm{d}Y_t,\quad t\in\mathbb{R} \tag{31}$$

with the control rate

$$c_t=-(\rho-u)X_t+\rho\mathrm{d}Y_t\,. \tag{32}$$

A mathematical difficulty arises in solving the BKE. For example, the infinite-dimensional analog of the BKE (20) is obtained by the following reformulation: with sufficiently smooth $y:[0,+\infty)\to\mathbb{R}$, such as $y\in C^1[0,+\infty)\cap L^2(0,+\infty)$ (continuously differentiability and square-integrability, both of which are necessary to well-define a quadratic solution guessed below), we can guess the formal convergence:

Discrete to continuous probability $c_i\to\pi(\mathrm{d}r)=p(r)\mathrm{d}r$ ($p$: density associated with $\pi$), (33)

Sequence to point values of a function $y_i\to y(r_i)$ ($i=1,2,3,\ldots$), (34)

Finite- to infinite-dimensional arguments $\Phi(x, y_1, \ldots, y_n) \to \Phi(x, y)$, $x \in \mathbb{R}$, $y \in L^2(0, +\infty)$, (35)

and taking the limit of $n \to +\infty$ in each term of the BKE, the following is obtained:

$$\left\{-(\rho - u)x + \sum_{i=1}^{n}(\rho - r_i)y_i\right\}\frac{\partial \Phi}{\partial x} \to \left\{-(\rho - u)x + \int_0^{+\infty}(\rho - \tau)y(\tau)\mathrm{d}\tau\right\}\frac{\partial \Phi(s, x, y(\cdot))}{\partial x}, \tag{36}$$

$$-\sum_{i=1}^{n} r_i y_i \frac{\partial \Phi}{\partial y_i} \to -\int_0^{+\infty} r y(r) \nabla_y \Phi(s, x, y(\cdot))(r)\mathrm{d}r, \tag{37}$$

$$\begin{aligned}&\sum_{i=1}^{n}(c_i A + B r_i y_i)\int_0^{+\infty}\left(\Phi(s, x+z, \ldots, y_i + z, \ldots) - \Phi\right)v(\mathrm{d}z)\\ &\to \int_0^{+\infty}\left(Ap(r) + Bry(r)\right)\int_0^{+\infty}\left(\Phi\left(s, x+z(r), \left(y + \delta(r - \cdot)z\right)(\cdot)\right) - \Phi(s, x, y(\cdot))\right)v(\mathrm{d}z(r))\mathrm{d}r\end{aligned}. \tag{38}$$

Here, $\nabla_y \Phi(s, x, y(\cdot))(r)$ is the Fréchet Derivative of $\Phi$ with respect to $y$ at $y(r)$, and $\delta$ is the Dirac's delta.

Consequently, we obtain the BKE of the infinite-dimensional problem

$$\begin{aligned}&-J + \left(x - \hat{X}\right)^2 + \left\{-(\rho - u)x + \int_0^{+\infty}(\rho - \tau)y(\tau)\mathrm{d}\tau\right\}\frac{\partial \Phi(s, x, y(\cdot))}{\partial x}\\ &-\int_0^{+\infty} r y(r)\nabla_y \Phi(s, x, y(\cdot))(r)\mathrm{d}r\\ &+\int_0^{+\infty}\left(Ap(r) + Bry(r)\right)\int_0^{+\infty}\begin{pmatrix}\Phi\left(s, x+z(r), \left(y + \delta(r - \cdot)z\right)(\cdot)\right)\\ -\Phi(s, x, y(\cdot))\end{pmatrix}v(\mathrm{d}z(r))\mathrm{d}r\end{aligned}, \quad x \in \mathbb{R}, \quad y \in L^2(0, +\infty). \tag{39}$$

This is an infinite-dimensional nonlocal PDE that needs to be carefully solved because of not only its infinite dimensions [41] but also the singular nonlocal term, the last term in (39) having Dirac's delta. A solution for this BKE is a pair of the constant $J \in \mathbb{R}$ and a smooth function $\Phi : \mathbb{R}^2 \times \left(C^1[0, +\infty) \cap L^2(0, +\infty)\right) \to \mathbb{R}$ such that all the terms of (39) are well defined.

Rigorous analysis of this BKE is nontrivial, but we can find a quadratic solution analogous to the finite-dimensional one. For any $x \in \mathbb{R}$ and $y \in C^1[0, +\infty) \cap L^2(0, +\infty)$, guess

$$\begin{aligned}\Phi(s, x, y) = {}&\frac{1}{2}Fx^2 + x\int_0^{+\infty} G(r)y(r)\mathrm{d}r + gx\\ &+\frac{1}{2}\int_0^{+\infty}\int_0^{+\infty}\Gamma(r, \tau)y(r)y(\tau)\mathrm{d}r\mathrm{d}\tau + \int_0^{+\infty}\gamma(r)y(\tau)\mathrm{d}r\end{aligned}. \tag{40}$$

Here, $F, g \in \mathbb{R}$ are constants, $G, \gamma : (0, +\infty) \to \mathbb{R}$ are smooth univariate functions, and $\Gamma : (0, +\infty)^2 \to \mathbb{R}$ is a smooth bivariate function. Then, we obtain

$$\nabla_y \Phi(x, y)(r) = xG(r)y(r) + \int_0^{+\infty}\Gamma(r, \tau)y(\tau)\mathrm{d}\tau + \gamma(r)y(r) \tag{41}$$

and

$$
\begin{aligned}
&\Phi\left(s, x+\delta(r-\cdot)z(\cdot), \left(y+\delta(r-\cdot)z\right)(\cdot)\right)-\Phi\left(s,x,y(\cdot)\right)\\
&=\frac{1}{2}F\left(x+z(r)\right)\left(x+z(r)\right)+\left(x+z(r)\right)\int_0^{+\infty}G(\tau)\left(y(\tau)+\delta(r-\tau)z(\tau)\right)\mathrm{d}w+g\left(x+z(r)\right)\\
&+\frac{1}{2}\int_0^{+\infty}\int_0^{+\infty}\Gamma(\tau,\omega)\left(y(\tau)+\delta(r-\tau)z(\tau)\right)\left(y(\omega)+\delta(r-\omega)z(\omega)\right)\mathrm{d}\tau\mathrm{d}\omega\\
&+\int_0^{+\infty}\gamma(\tau)\left(y(\tau)+\delta(r-\tau)z(\tau)\right)\mathrm{d}\tau\\
&-\left(\frac{1}{2}Fx^2+x\int_0^{+\infty}G(\tau)y(\tau)\mathrm{d}\tau+gx+\frac{1}{2}\int_0^{+\infty}\int_0^{+\infty}\Gamma(\tau,\omega)y(\tau)y(\omega)\mathrm{d}\tau\mathrm{d}\omega+\int_0^{+\infty}\gamma(\tau)y(\tau)\mathrm{d}\tau\right)\\
&=\left\{z(r)\right\}^2\left\{\frac{1}{2}F+G(r)+\frac{1}{2}\left\{z(r)\right\}^2\right\}+z(r)\left\{g+\gamma(r)\right\}\\
&+z(r)x\left\{F+G(r)\right\}+z(r)\int_0^{+\infty}\left\{G(\tau)+\Gamma(\tau,r)\right\}y(\tau)\mathrm{d}\tau
\end{aligned}
\quad . \tag{42}
$$

These are the infinite-dimensional counterparts of (A.1) and (A.2) in **Appendix A.1**. The other parts can be evaluated similarly. Technically, we find the coefficients $F, g, G, \gamma, \Gamma$ and the correspondence between $F$ and $a_{0,0}$, $G(r_i)$ and $a_{i,0}$ ($1 \le i \le n$), $\Gamma(r_i, r_j)$ and $a_{i,j}$ ($1 \le i, j \le n$), $g$ and $b_0$, and $\gamma(r_i)$ and $b_i$ ($1 \le i \le n$). The smooth function (40) is, therefore, a point-wise solution for the BKE in the infinite-dimensional case. Through a simple calculation analogous to **Appendix A.1** with each summation read as an appropriate integral with respect to $\pi(\mathrm{d}r)$ (i.e., $\sum_{i=1}^{n} c_i f(r_i) \to g = \int_0^{+\infty} f(r)\pi(\mathrm{d}r)$), we obtain

$$
F=\frac{1}{\rho-u},\quad G(r)=\frac{1}{\rho-u}\frac{\rho-rD}{rD+\rho-u},\quad \Gamma(r,\tau)=\frac{(\rho-rD)(\rho-\tau D)}{(r+\tau)D}\left\{\frac{1}{rD+\rho-u}+\frac{1}{\tau D+\rho-u}\right\}, \tag{43}
$$

$$
g=-\frac{2}{\rho-u}\hat{X}+\frac{1}{\rho-u}(q+1)AM_1\int_0^{+\infty}\frac{1}{r}\frac{r}{rD+\rho-u}\pi(\mathrm{d}r), \tag{44}
$$

$$
\begin{aligned}
&r\left(1-BM_1\right)\left(\gamma(r)+g\right)\\
&=-2q\hat{X}+q(q+1)AM_1\int_0^{+\infty}\frac{1}{r}\frac{r}{rD+\rho-u}\pi(\mathrm{d}r)+\frac{1}{2(\rho-u)}BM_2\frac{(\rho-u)rD+\rho^2}{D(rD+\rho-u)}\\
&+\frac{1}{\rho-u}AM_1\frac{\rho-rD}{rD+\rho-u}+\frac{1}{\rho-u}AM_1\int_0^{+\infty}\frac{(\rho-rD)(\rho-\tau D)}{(r+\tau)D}\left\{\frac{1}{rD+\rho-u}+\frac{1}{\tau D+\rho-u}\right\}\pi(\mathrm{d}r)
\end{aligned}
\quad , \tag{45}
$$

and consequently

$$
J=\left(\hat{X}-q\mathbb{E}[Y]\right)^2+\mathrm{Var}[Y]\left(\int_0^{+\infty}\frac{1}{r}\pi(\mathrm{d}r)\right)^{-1}\int_0^{+\infty}\frac{1}{r}\frac{r+q^2(\rho-u)}{r+\rho-u}\pi(\mathrm{d}r), \tag{46}
$$

where $\mathbb{E}[Y]=\frac{AM_1}{D}\int_0^{+\infty}\frac{1}{r}\pi(\mathrm{d}r)$ and $\mathrm{Var}[Y]=\frac{AM_2}{2D}\int_0^{+\infty}\frac{1}{r}\pi(\mathrm{d}r)$.

Notably, (46) is obtained by formally letting $n \to +\infty$ in (21). Convergence is true if the sequences $\{c_i\}_{1\le i\le n}$ and $\{r_i\}_{1\le i\le n}$ satisfy some conditions. We present an example in the next subsection, and here, we assume that such $\{c_i\}_{1\le i\le n}$ and $\{r_i\}_{1\le i\le n}$ are specified. Then, under $n \to +\infty$ we obtain

$$K=(\rho-u)^2(1-q)^2\,\mathrm{Var}[Y]\left(\int_0^{+\infty}\frac{1}{r}\pi(\mathrm{d}r)\right)^{-1}\int_0^{+\infty}\frac{D}{rD+\rho-u}\pi(\mathrm{d}r). \tag{47}$$

Consequently, we obtain the following solution to the minimization problem of $J$ under $q=\hat{X}/\mathbb{E}[Y]$.

***Proposition 8:*** *Assume $q=\hat{X}/\mathbb{E}[Y]$ and the dynamics (30)–(32). Then, the minimization problem of $J$ subject to (16) and (17) under $q=\hat{X}/\mathbb{E}[Y]$ is solved as follows:*

*(d) If $q=1$, then $\min J=\mathrm{Var}[Y]$. The minimizing couple $(\rho,u)$ is $u=0$ and arbitrary $\rho>0$.*

*(e) If $q>1$ **(Water Adding case)**, then $\inf_{h=\rho-u} J(h)=J(0)=\mathrm{Var}[Y]$. The minimizing couple does not exist, but $\lim_{h\to+0} J(h)=\mathrm{Var}[Y]$.*

*(f) If $0<q<1$ **(Water Abstracting case)**, then $\inf_{h=\rho-u} J(h)=J(\bar{h})=\mathrm{Var}[Y]$ with $\bar{h}$ uniquely determined from (48). The minimizing couple $(\rho,u)$ satisfies $q=\dfrac{\rho}{\rho-u}$ and $\rho-u=\bar{h}$, and the minimum of $J$ is achieved by $(\rho,u)=(q\bar{h},-(1-q)\bar{h})$:*

$$K(\bar{h})=(1-q)^2\,\mathrm{Var}[Y]\left(\int_0^{+\infty}\frac{1}{r}\pi(\mathrm{d}r)\right)^{-1}\bar{h}^2\int_0^{+\infty}\frac{D}{rD+\bar{h}}\pi(\mathrm{d}r)=\bar{K}. \tag{48}$$

### 3.5 Approximation scheme

We used a sufficiently fine Markovian lift to approximate the integrals $J,K$ and the nonlinear equation (48). The problem in this step is how to numerically approximate integrals that appear in the representations of $J,K$. Here, we present a convergent discretization.

First, the sequences $\{c_i\}_{1\le i\le n}$ and $\{r_i\}_{1\le i\le n}$, are chosen so that the finite-dimensional supCBI process converges to the original supCBI process, at least considering the distribution. According to [32], this convergence holds if

$$\left|\sum_{i=1}^{n}\frac{c_i}{r_i}-\int_0^{+\infty}\frac{1}{r}\pi(\mathrm{d}r)\right|\to+0 \quad\text{as}\quad n\to+\infty. \tag{49}$$

Herein, we constructed such a discretization method based on quantiles, and it suffices in approximating the above-mentioned integrals. Without significant loss of generality, we assume that $\pi$ admits a smooth and absolutely continuous density function and satisfies $\int_0^{+\infty}\frac{1}{r}\pi(\mathrm{d}r)<+\infty$. These are satisfied in our application.

We present a discretization method without truncating the domain. First, we set the discrete probability measure $\pi_n(\mathrm{d}r)=\sum_{i=1}^{n}c_i\delta(r-r_i)$. Given $n=2^m$ ($m\in\mathbb{N}$) and $k=0,1,2,\ldots,n$, we set the quantile $\theta_{n,k}$ ($1 \le k \le n$) uniquely determined from

$$\frac{k}{n} = \int_0^{\theta_{n,k}} \pi(\mathrm{d}r). \tag{50}$$

Notably, $\theta_{n,n} = +\infty$. Then, we set $r_i = \theta_{2n,2i-1}$ and $c_i = \int_{\theta_{2n,2i-2}}^{\theta_{2n,2i}} \pi(\mathrm{d}r)$ ($1 \le i \le n$). From (50), we obtain

$$c_i = \int_0^{\theta_{2n,2i}} \pi(\mathrm{d}r) - \int_0^{\theta_{2n,2i-2}} \pi(\mathrm{d}r) = \frac{2i}{2n} - \frac{2i-2}{2n} = \frac{1}{n} \quad (1 \le i \le n). \tag{51}$$

The sequence $\{c_i\}_{1 \le i \le n}$ in this case is constant. In addition, for $r \in \left(\theta_{2n,2i-21}, \theta_{2n,2i}\right)$, we have

$$\left| \int_0^r \pi(\mathrm{d}r) - \int_0^r \pi_n(\mathrm{d}r) \right| \le \max\left\{ \left| \frac{2i}{2n} - r \right|, \left| \frac{2i-2}{2n} - r \right| \right\} \le \left( \frac{2i}{2n} - \frac{2i-2}{2n} \right) = \frac{1}{n}. \tag{52}$$

It follows that the cumulative function $\int_0^{\bullet} \pi_n(\mathrm{d}r)$ uniformly converges to $\int_0^{\bullet} \pi(\mathrm{d}r)$ in $(0,+\infty)$. As a byproduct, by the Portemanteau theorem (Theorem 13.16 of Klenke [42]), we obtain (Recall $n = 2^m$)

$$\left| \int_0^{+\infty} f(r)\pi(\mathrm{d}r) - \int_0^{+\infty} f(r)\pi_n(\mathrm{d}r) \right| \to 0 \quad \text{as} \quad m \to +\infty \tag{53}$$

for any bounded and continuous functions $f$ on $[0,+\infty)$.

Then, for each constant $\kappa > 0$, we obtain

$$\left| \int_0^{+\infty} \frac{1}{\max\{r,\kappa\}} \pi(\mathrm{d}r) - \int_0^{+\infty} \frac{1}{\max\{r,\kappa\}} \pi_n(\mathrm{d}r) \right| \to 0 \quad \text{as} \quad m \to +\infty \tag{54}$$

because the function $\frac{1}{\max\{r,\kappa\}}$ ($r > 0$) is nonnegative and bounded from above by $\kappa^{-1}$. In addition,

$$\left| \int_0^{+\infty} \frac{1}{r} \pi(\mathrm{d}r) - \int_0^{+\infty} \frac{1}{\max\{r,\kappa\}} \pi(\mathrm{d}r) \right| \to 0 \quad \text{as} \quad \kappa \to +0 \tag{55}$$

by $\int_0^{+\infty} \frac{1}{r} \pi(\mathrm{d}r) < +\infty$. We have

$$\begin{aligned} &\left| \int_0^{+\infty} \frac{1}{\max\{r,\kappa\}} \pi_n(\mathrm{d}r) - \int_0^{+\infty} \frac{1}{r} \pi(\mathrm{d}r) \right| \\ &\le \left| \int_0^{+\infty} \frac{1}{\max\{r,\kappa\}} \pi_n(\mathrm{d}r) - \int_0^{+\infty} \frac{1}{\max\{r,\kappa\}} \pi(\mathrm{d}r) \right| + \left| \int_0^{+\infty} \frac{1}{\max\{r,\kappa\}} \pi(\mathrm{d}r) - \int_0^{+\infty} \frac{1}{r} \pi(\mathrm{d}r) \right|. \end{aligned} \tag{56}$$

Then, we obtain

$$\limsup_{\kappa \to +0} \limsup_{m \to +\infty} \left| \int_0^{+\infty} \frac{1}{\max\{r,\kappa\}} \pi_n(\mathrm{d}r) - \int_0^{+\infty} \frac{1}{r} \pi(\mathrm{d}r) \right| = 0, \tag{57}$$

showing

$$\lim_{\kappa \to +0} \lim_{m \to +\infty} \int_0^{+\infty} \frac{1}{\max\{r,\kappa\}} \pi_n(\mathrm{d}r) = \int_0^{+\infty} \frac{1}{r} \pi(\mathrm{d}r). \tag{58}$$

We choose a sequence $\{\kappa_m\}_{m=1,2,3,\ldots}$ converging to 0 such that $0 < \kappa_m < \theta_{2^{m+1},1}$. Substituting $\kappa_m$ for $\kappa$,

we obtain

$$\int_0^{+\infty} \frac{1}{\max\{r,\kappa_m\}} \pi_n(\mathrm{d}r) = \int_0^{+\infty} \sum_{i=1}^{n} \frac{1}{\max\{r_i,\kappa_m\}} \delta(r-r_i)\mathrm{d}r = \int_0^{+\infty} \sum_{i=1}^{n} \frac{1}{r_i} \delta(r-r_i)\mathrm{d}r = \int_0^{+\infty} \frac{1}{r} \pi_n(\mathrm{d}r). \quad (59)$$

From the convergence result (57)–(58) and the equality (59), we obtain

$$\lim_{m\to+\infty} \int_0^{+\infty} \frac{1}{\max\{r,\kappa_m\}} \pi_n(\mathrm{d}r) = \lim_{m\to+\infty} \int_0^{+\infty} \frac{1}{r} \pi_n(\mathrm{d}r) = \int_0^{+\infty} \frac{1}{r} \pi(\mathrm{d}r) \quad (60)$$

since $m = \log_2 n$. Notably, the sequence $\{\kappa_m\}_{m=1,2,3,...}$ was introduced only for a technical reason, and it does not appear in the implementation of the discretization.

Consequently, the above-specified sequences $\{c_i\}_{1\le i\le n}$ and $\{r_i\}_{1\le i\le n}$ suffice for our purpose. In our application, we use a Gamma type $\pi$. Each quantile $\theta_{n,k}$ is numerically computed with an arbitrary level of accuracy using the sequential representation of incomplete gamma functions (Lemma 14 of Greengard and Rokhlin [43]) combined with a bisection method to solve the nonlinear equation (50) for $\theta_{n,k}$. The convergence of the integrals $J, K$ is similar as their integrands are strictly bounded or have a singularity similar to $r^{-1}$. The computational results in **Appendix B** are consistent with the theoretical result: the integral $R = \int_0^{+\infty} r^{-1}\pi(\mathrm{d}r)$ can be approximated by the summation $R_n = \sum_{i=1}^{n} c_i r_i^{-1}$.

### 3.6 Solution algorithm

This sub-section explains a solution algorithm for numerically solving our optimization problem, which is (15) subject to constraints. The supCBI process is firstly identified using available data. Secondly, the probability measure $\pi$ is discretized using the Markovian lift explained in **Section 3.5**. Thirdly, set parameter values of the objective and solve (29) with respect to $\bar{h}$ using a common method for nonlinear equations such as a Picard method. Finally, we can find optimal $(\rho, u)$ from **Proposition 7**. Our optimization problem is therefore solvable using a common method and does not require the use of sophisticated methodologies. This practical implementability is an advantage of our optimization problem.

## 4. Applications

We identified the supCBI process at a midstream reach of an existing river and solved the optimization problem of the Water Abstraction case as it offers feasible solutions of interest.

### 4.1 Study site

We applied the optimization problem to an existing river environment. The study site is a midstream reach of the Hii River flowing in the eastern part of Shimane Prefecture in Japan (**Figure 2**). We have studied the river as a representative river environment in a mountainous region where balancing between water resources for human lives (e.g., hydropower generation and irrigation) and environmentally required streamflows for maintaining healthy aquatic ecosystems, including several fishery resources, has recently been a huge concern [44].

The river reach studied in this herein is 16 (km) long and is bounded by the Obara Multipurpose Dam (O Dam, managed by the Ministry of Land, Infrastructure, Transport and Tourism (MLIT)) at the upstream end and Hinobori Weir (H Weir) for trapping soil materials at the downstream end. This reach is a famous recreational fishing area of *Plecoglossus altivelis*, a recreationally important iconic fishery resource in Japan [45]. Therefore, very large river discharge is dangerous for anglers, and small discharge results in the lack of the shear stress acting on the riverbed and causes the bloom of nuisance periphyton [46], which affects the ecosystem balance and fishing utilities. In the study area, the recent growing need for a cleaner energy supply demands the construction of another hydropower station to complement the one upstream of O Dam. Water abstraction for irrigation would also increase in the future due to the need for clean electricity supplies in the surrounding watershed area. Therefore, streamflow management at the study site is a crucial issue.

We chose this study site also because the hydropower station of a private enterprise abstracting the river water upstream of O Dam has been idling since the beginning of June 2020 [47]. The river water flowing into the O Dam and, hence, the water discharged from the dam has increased since the beginning of June 2020, resulting in an artificial regime shift in the study site. Thus, we could consider two streamflow conditions before and after the operational change in the hydropower station. More specifically, the average discharge downstream of the O Dam increased at the beginning of June 2020.

The public data of hourly outflow discharge of O Dam since April 1, 2016, are available (MLIT, http://www1.river.go.jp/cgi-bin/SrchDamData.exe?ID=607041287705020&KIND=1&PAGE=0). Thus, we considered the two periods, which are the first four years (Period 1: June 1, 2016, to May 31, 2020) and the second two years (Period 2: June 1, 2020, to May 31, 2022). Because we are interested in the 16-km reach, which cannot be considered a single point, we employed hourly discharge time series at each point in the reach using a longitudinally one-dimensional (1D) hydrodynamic model, as explained below.

### 4.2 Flow field

The discharge time series at each point in the target reach was emulated by the hydrodynamic model driven using the record of the outflow discharge of the O Dam as the inflow boundary conditions at the upstream end, a weir boundary condition at the downstream end, and lateral inflows based on a distributed hydrological model [48]. This model has been verified using water depth data collected in 2021, but we reexamined it herein using the latest water depth time-series data collected at the study site at a 10-min interval using a water-level gauge (HOBO, U-20). The 1D computational grid contains 302 grid points (Points 1 through 302 from the upstream), and the grid interval to discretize the river reach is approximately 50 m with the time step so that the Courant–Friedrichs–Lewy number is 0.5 at each time step. Considering the inflow boundary condition and lateral inflows as stochastic processes, emulating the flow field in the river reach results in numerically solving a stochastic PDE in a 1D domain. **Figure 3** compares the measured and computed water depth of Station Y at 10.2 (km) downstream of O Dam, and they agree well during both high- and low-flow periods. This comparison suggests that the emulated flow field can be used to identify the MMA model in this river reach.

We statistically evaluated the flow field in Periods 1 and 2. **Figures 4(a)–(d)** compare the average ($m^3/s$), variance ($m^6/s^2$), skewness, and kurtosis of the discharge at each computational grid point in the domain in Period 1, assuming stationarity. Similarly, **Figures 5(a)–(d)** compare the statistics in Period 2, assuming stationarity. The average discharge in Period 2 is larger than that of Period 1; the difference is approximately 7 ($m^3/s$), indicating more than a 75% increase from Period 1. The increase in discharge in Period 2 is consistent with the change in the inflow discharge to O Dam [47]. This supports our hypothesis that the average discharge downstream of the dam is larger in Period 2 than it is in Period 1. We also analyzed other statistics. The increase in variance in Period 1 from Period 2 is 30%–60%, and those of skewness and kurtosis are about −8% to −20% and −16% to −34%, respectively. This implies that idling the upstream hydropower station shifted the flow field to a statistically more variable, less skewed, and less sharp one.

On the spatial distributions of the statistics, **Figures 4 and 5** show that the target reach can be divided into three parts, where the statistics vary gradually (1.2–3.0 and 5.1–16 (km) from the dam) and more rapidly (3.0–5.1 (km) from the dam). The more rapid variation was considered due to the larger volume of lateral inflows received in the middle reach at the study site. Considering this computational result, in the application of the proposed model below, it is not necessary to consider all the grid points but only some of them with a sparser interval of approximately 1 (km) (**Figure 6;** these points are hereafter specific grid points for simplicity). Finally, in both periods, the ACF at each point has an algebraic decay as predicted from the MMA process, as discussed in more detail in the next subsection.

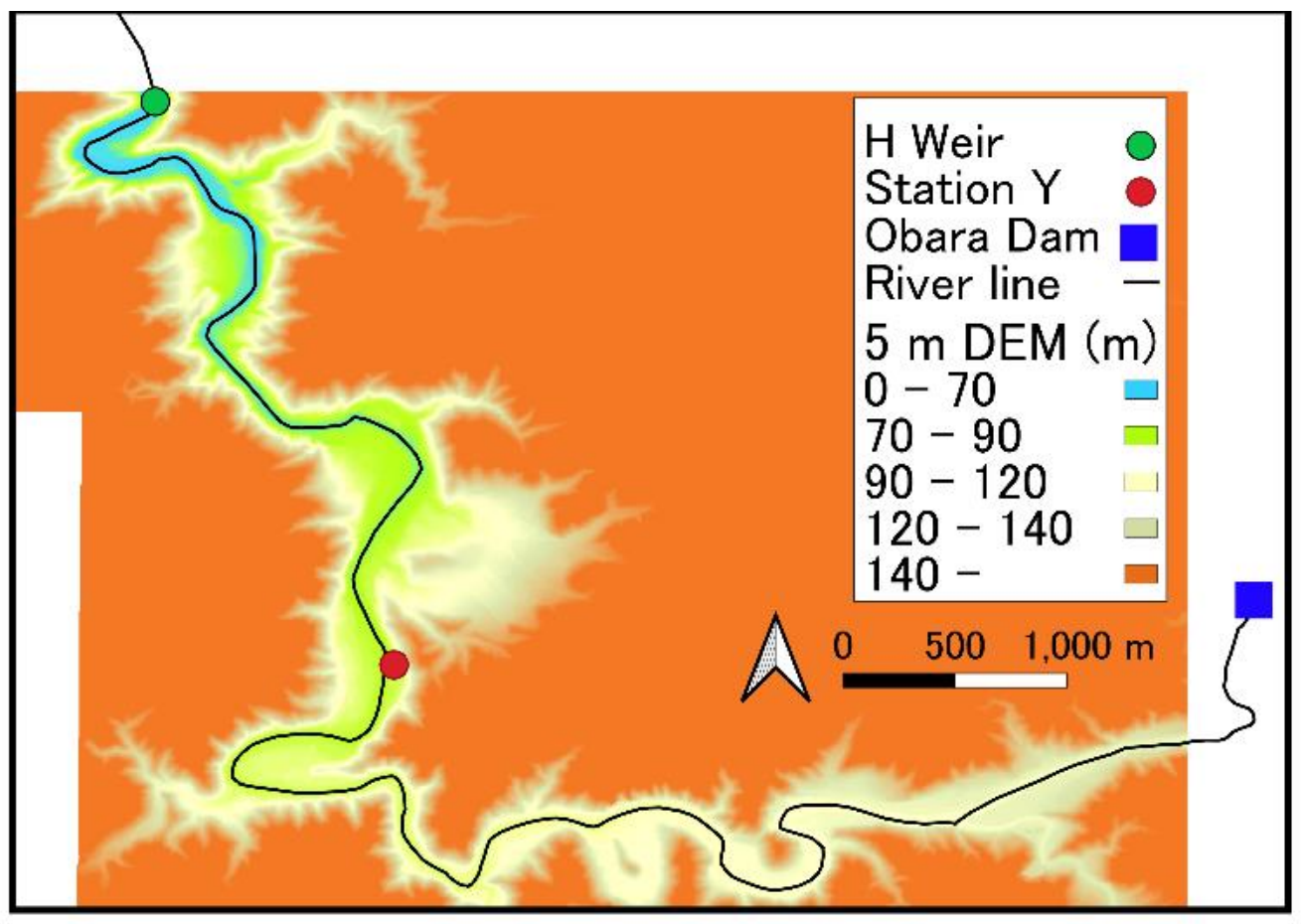

**Figure 2.** Map of the site studied in this study. The elevation data was constructed using a 5-m DEM provided by the Geospatial Information Authority of Japan (https://fgd.gsi.go.jp/download/mapGis.php?tab=dem).

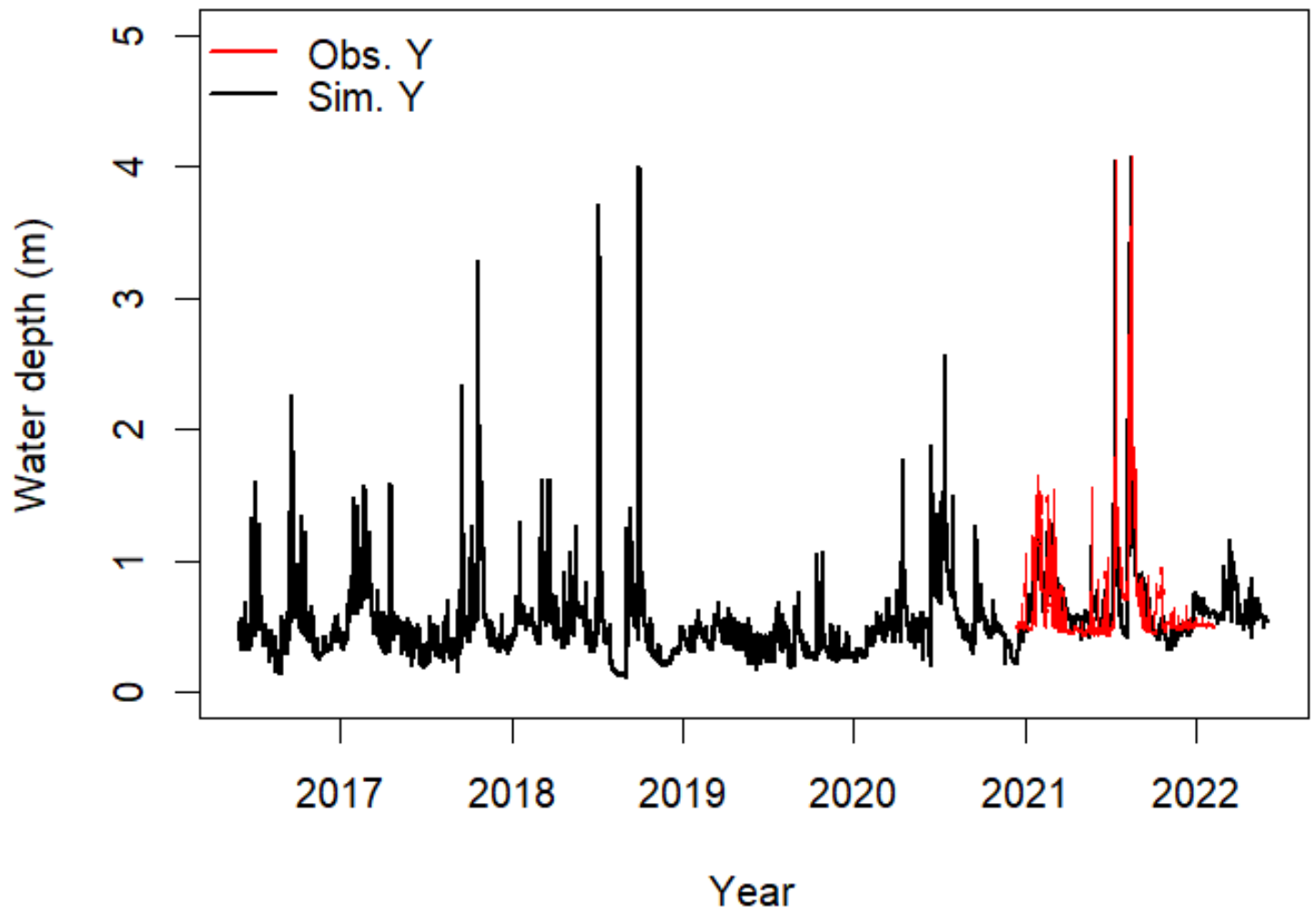

**Figure 3.** Comparison of measured and computed water depth at Station Y at 10.2 (km) downstream of O Dam.

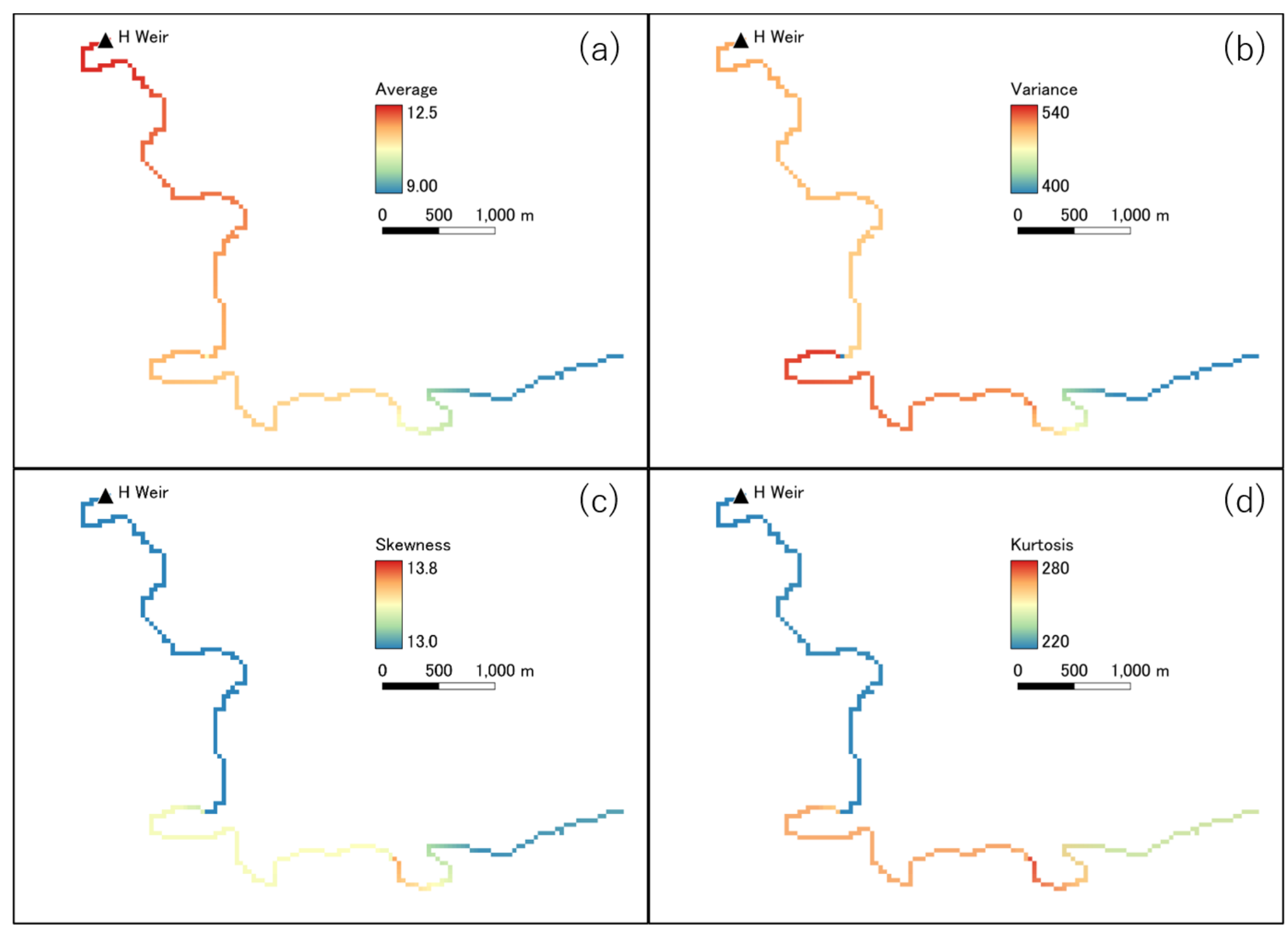

**Figure 4.** Computed stationary (a) average ($m^3/s$), (b) variance ($m^6/s^2$), (c) skewness (-), and (d) kurtosis (-) in the river reach during Period 1.

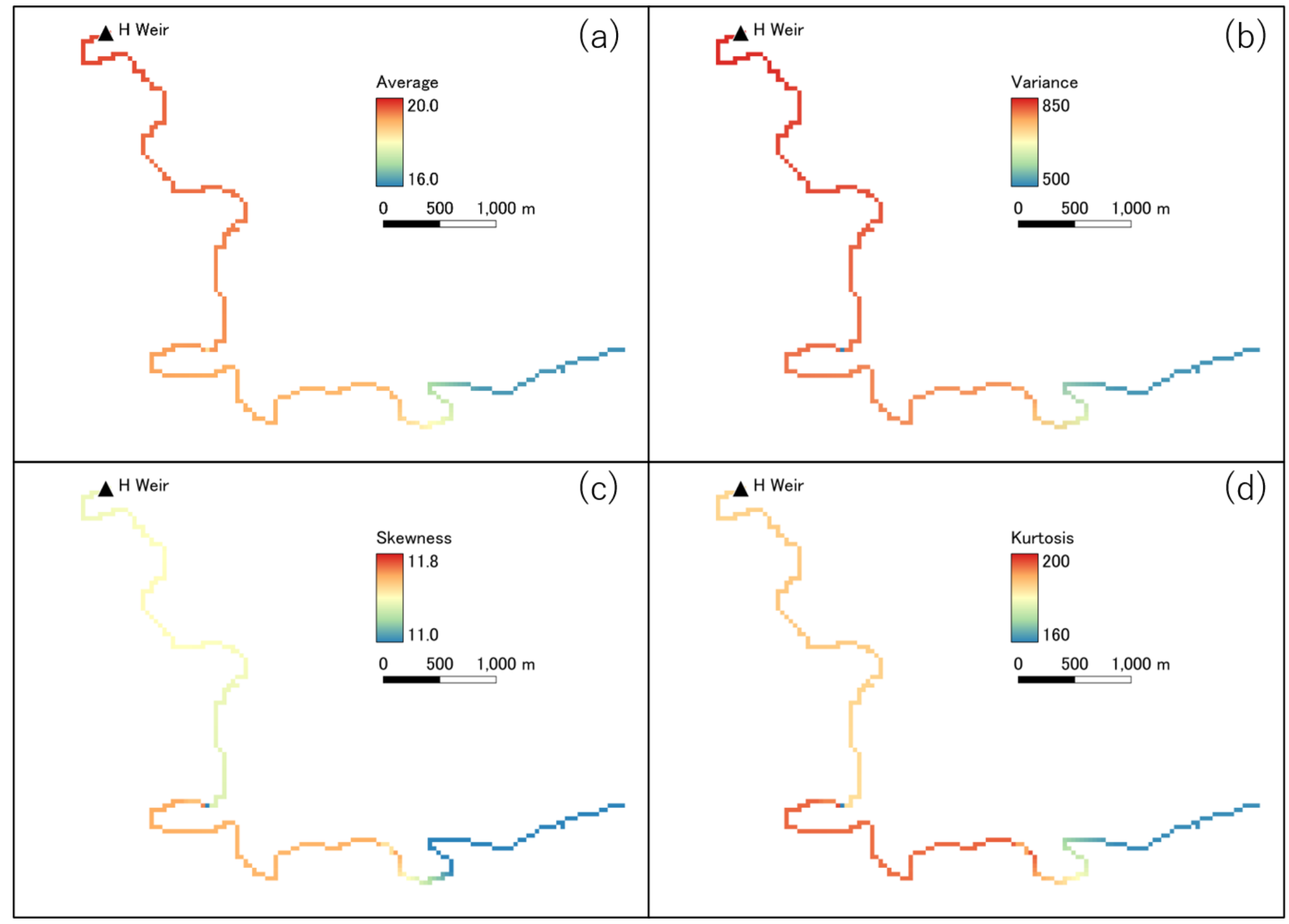

**Figure 5.** Computed stationary (a) average ($m^3/s$), (b) variance ($m^6/s^2$), (c) skewness (-), and (d) kurtosis (-) in the river reach during Period 2.

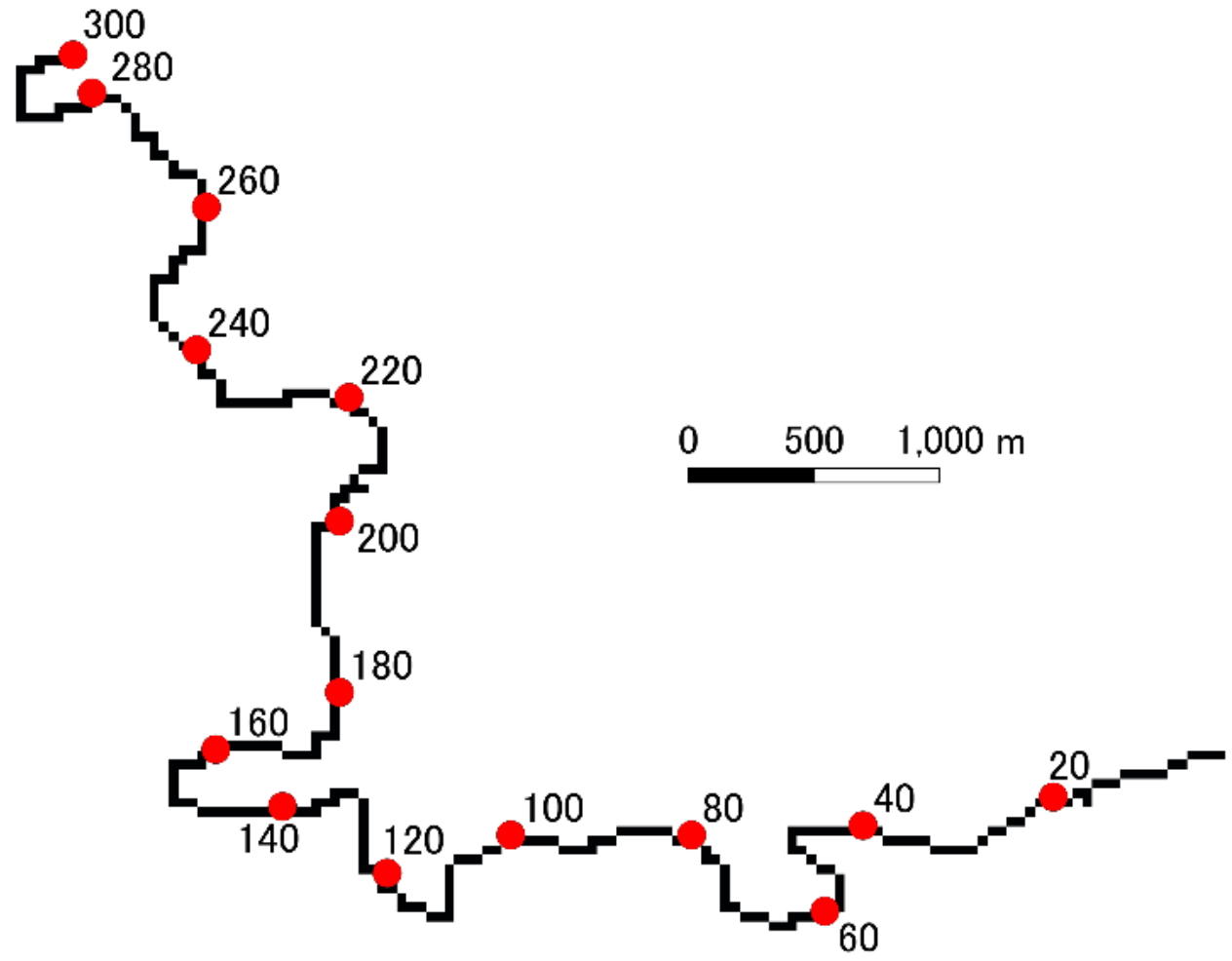

**Figure 6.** Locations of the chosen sparser points: Points 20 through 300.

### 4.3 Model identification

The MMA process is identified at each grid point following the procedure reported by [17, 32] with a slight modification. First, for a fixed $D \in (0,1)$, the probability measure $\pi$, including the parameters $\alpha$ and $\beta$, is identified using (6) and a common least-square method. In this step, the data with the minimal delay time interval containing 0 such that the ACF remains positive is used. Second, the Lévy measure $\nu$, namely, the remaining parameters, is identified. This is conducted using the least-square method so that the following error metric is minimized:

$$\begin{aligned} E = &\left(\frac{\mathrm{Average}_{\mathrm{M}} - \mathrm{Average}_{\mathrm{D}}}{\mathrm{Average}_{\mathrm{D}}}\right)^2 + \left(\frac{\mathrm{Variance}_{\mathrm{M}} - \mathrm{Variance}_{\mathrm{D}}}{\mathrm{Variance}_{\mathrm{D}}}\right)^2 \\ &+ \left(\frac{\mathrm{Skewness}_{\mathrm{M}} - \mathrm{Skewness}_{\mathrm{D}}}{\mathrm{Skewness}_{\mathrm{D}}}\right)^2 + \left(\frac{\mathrm{Kurtosis}_{\mathrm{M}} - \mathrm{Kurtosis}_{\mathrm{D}}}{\mathrm{Kurtosis}_{\mathrm{D}}}\right)^2 . \end{aligned} \quad (61)$$

Subscripts "M" and "D" indicate the model and data, respectively. In our previous studies, we used standard deviation instead of variance. Herein, we use variance as it is a cumulant, whereas the standard deviation used in [32] is not, implying that it is more natural to use variance. In addition, variance is directly related to the objective of our optimization problem because it can be identified as the objective of it. Notably, the model identification results presented below are almost the same if we choose the standard deviation in the error metric. Here, we set $D = 0.5$. The obtained results are not significantly affected by the choice of $D$. We assume a tempered stable model $\nu(\mathrm{d}z) = \exp(-c_2 z) z^{-(1+c_1)} \mathrm{d}z$ ($c_2 > 0$, $c_1 < 1$) as the simplest model to describe subordinators having bounded and unbounded variations in a unified manner. The probability measure $\pi$ is assumed to be a Gamma type $\pi(\mathrm{d}r) \sim r^{\alpha-1} e^{-r/\beta} \mathrm{d}r$ introduced earlier due to its simplicity and capability of generating long-memory MMA processes.

The identification results are summarized in **Appendix C**. **Table C.1** lists the identification results for Period 1 at three points in the target reach, which are in the upstream region where the statistics vary gradually (Point 20, 2.3 (km) from the dam), the region where the statistics rapidly vary (Point 60, 4.3 (km) from the dam), and the downstream region where the statistics gradually vary (Point 180, 10.2 (km) from the dam). Similarly, **Table C.2** lists the identification results in Period 2. For Point 20, the ACF in Period 2 is a long-memory type ($\alpha \leq 2$), whereas that of Period 1 is not ($\alpha > 2$), implying that the regime shift in the target reach affects the temporal correlation structure (**Figures C.1** and **C.2** in **Appendix C**). Similar results were obtained for other points and are not reported here.

In summary, the MMA process was successfully fitted into the target reach. The statistical analysis combined with the hydraulic simulation suggests that, in each period, the statistics (average, variance, skewness, and kurtosis) depend on the distance from the dam, whereas the autocorrelation structure (ACF) of the discharge does not. Therefore, the memory effect was inherited along the target reach in each period.

### 4.4 Computation and analysis

The optimization problem was solved at each point in the target river reach. We focused on the Abstracting Water case as it is more realistic than the Water Adding case. We used the technique highlighted in **Remark 3** and focused on the abstraction $Q_{\mathrm{abs}}$ rather than the target discharge $\hat{Q}$. Thus, we used the following one-to-one relationship between $\hat{Q}$ and $\underline{Q}$:

$$\hat{X} = \hat{Q} - q\underline{Q} = (1-q)\underline{Q} + \mathbb{E}[Y_n] - Q_{\mathrm{abs}} \quad \text{or equivalently} \quad Q_{\mathrm{abs}} = (1-q)\left(\underline{Q} + \mathbb{E}[Y_n]\right). \tag{62}$$

Before performing a detailed analysis, we checked the convergence of the discretization method based on the Markovian lift. We used the degree of freedom $n = 2^{13} = 8192$, with which sufficiently good results were obtained, as suggested by the numerical experiments in **Appendix B**.

First, we focused on the problem at a single point, which is Point 180, located 10.2 km downstream from O Dam. **Figures 7–9** show the computed minimized variance $J$ and the corresponding optimal $\rho$ and $u$ with respect to the cost constraint $\bar{K}$. Here, we considered different values of the abstraction $Q_{\mathrm{abs}}$ (m$^3$/s) (3, 5, and 7 (m$^3$/s)), which are realistic because the minimum flow required in the target reach is 1.0 (m$^3$/s) considering the average in **Tables C.1 and C.2**. **Figure 7** shows that, for each case and period, the $\bar{K}$-$J$ curve is a continuum of the minimized variance given the cost constraint. The curve has a convex shape and approaches the theoretical lower bound expressed in (28) as the cost constraint decreases ($\bar{K}$ increases). As shown in **Figures 9 and 10**, the dependence of optimal $\rho$ (resp., $u$) on $\bar{K}$ increases almost linearly (resp., concave and almost linearly decreasing). A sharp decrease in the minimized variance is expected only when the cost constraint is small with these optimal parameters. The dependence of the minimized variance $J$ between Periods 1 and 2 is clear; for the same $\left(\bar{K}, Q_{\mathrm{abs}}\right)$, the minimized variance $J$ is smaller for Period 1. This is attributed to the larger variance of the discharge in Period 2 than in Period 1, indicating that the flow condition in Period 2 is more difficult to regulate around a target value.

We extended the analysis to each point in the target reach to show that the proposed model can be employed in creating water infrastructure to abstract river water. This was conducted by exploring the point where the minimized variance $J$ becomes minimum among all the candidate points in the target reach. We solved this nested optimization problem in each period for different values of the cost constraint $\bar{K}$. Similarly, we explored the maximizing point of the minimum variance $J$ where the worst-case performance is achieved among the points examined, which is a worst-case minimization problem.

**Figure 11** shows a plot of the computed minimum variance $J$ with respect to the cost constraint $\bar{K}$ for all the specific grid points. The target abstraction is fixed to $Q_{\mathrm{abs}} = 5$ ($\mathrm{m}^3$/s). **Table 1** lists the optimization results. The minimizing location for all examined $\bar{K}$ (0.001, 0.002, 0.003,…, 10.000 ($\mathrm{m}^6/\mathrm{s}^2$)) and $Q_{\mathrm{abs}}$ (1 through 8 ($\mathrm{m}^3$/s)) is the most upstream point (Point 20) where the variance becomes minimum in the target reach. In contrast, the worst point is Point 160 or 300, depending nonmonotonically on $\bar{K}$ and $Q_{\mathrm{abs}}$ (**Table 1**).

The nested optimization problem here does not consider costs other than the control cost, including the construction and maintenance costs of the infrastructure and water conveying channels. However, the proposed methodology can be incorporated into such an engineering analysis as a pivotal optimization model. Based on the proposed approach, we shall conduct spatiotemporal analysis to optimize both the water resources supply and environments in the study area considering its land use.

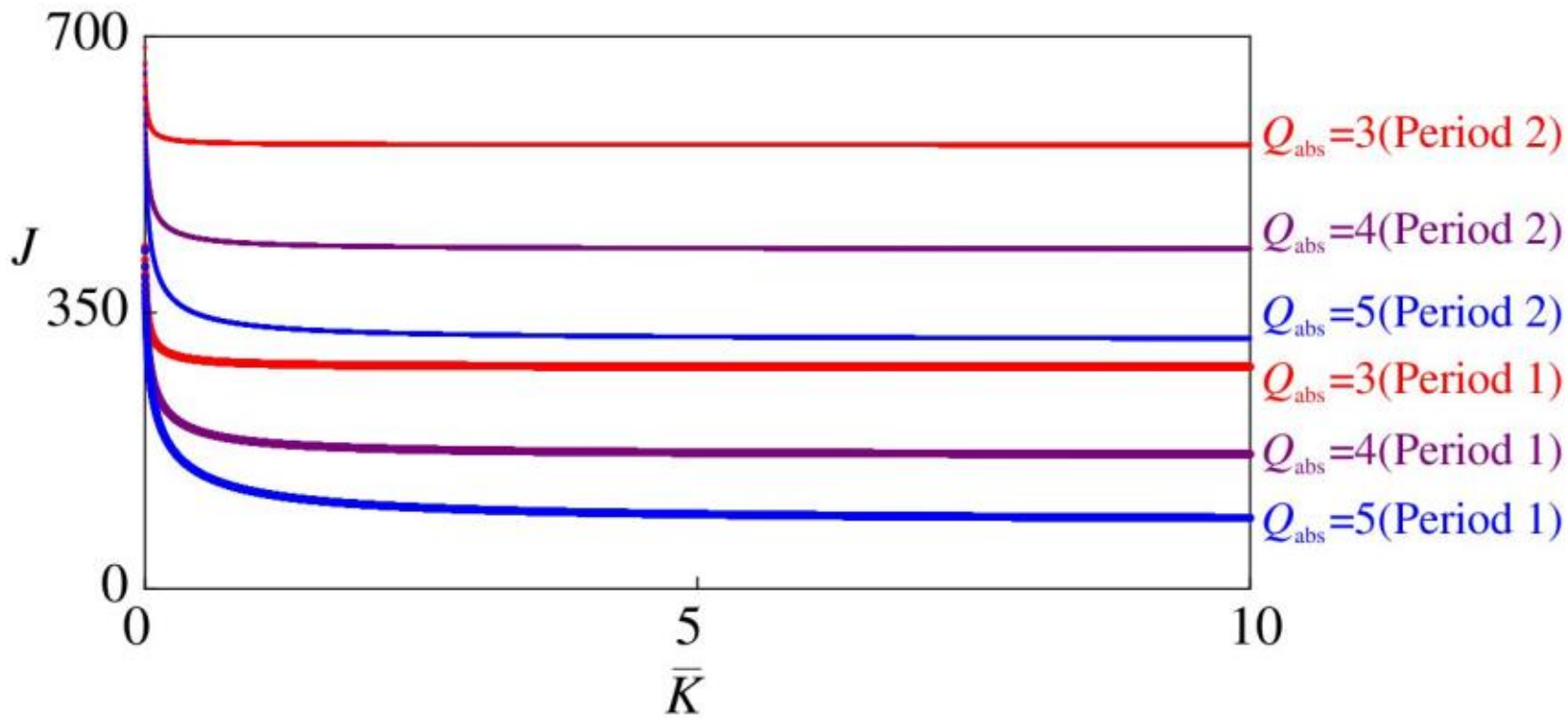

**Figure 7.** The computed minimum variance $J$ with respect to the cost constraint $\overline{K}$.

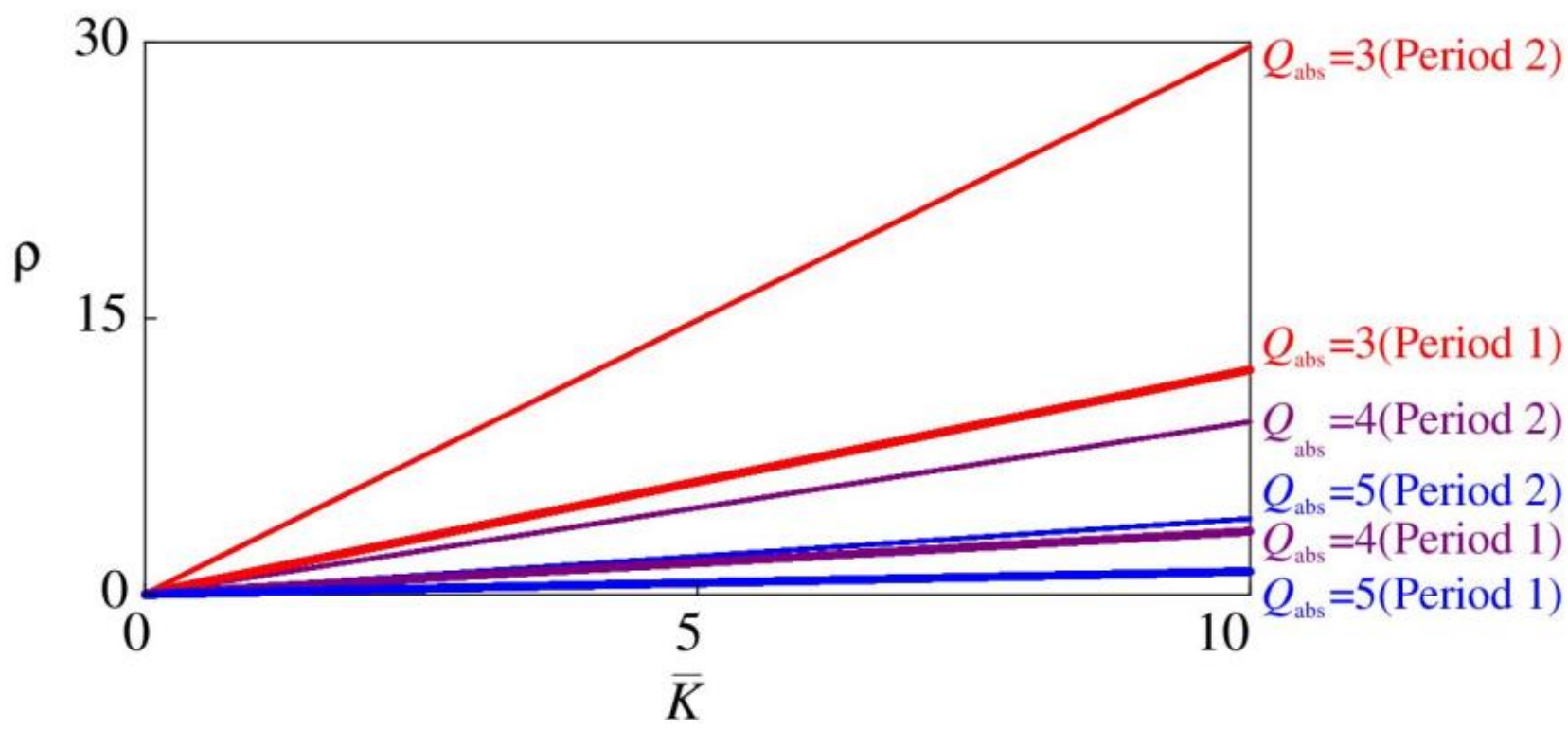

**Figure 8.** The computed optimal $\rho$ with respect to the cost constraint $\overline{K}$.

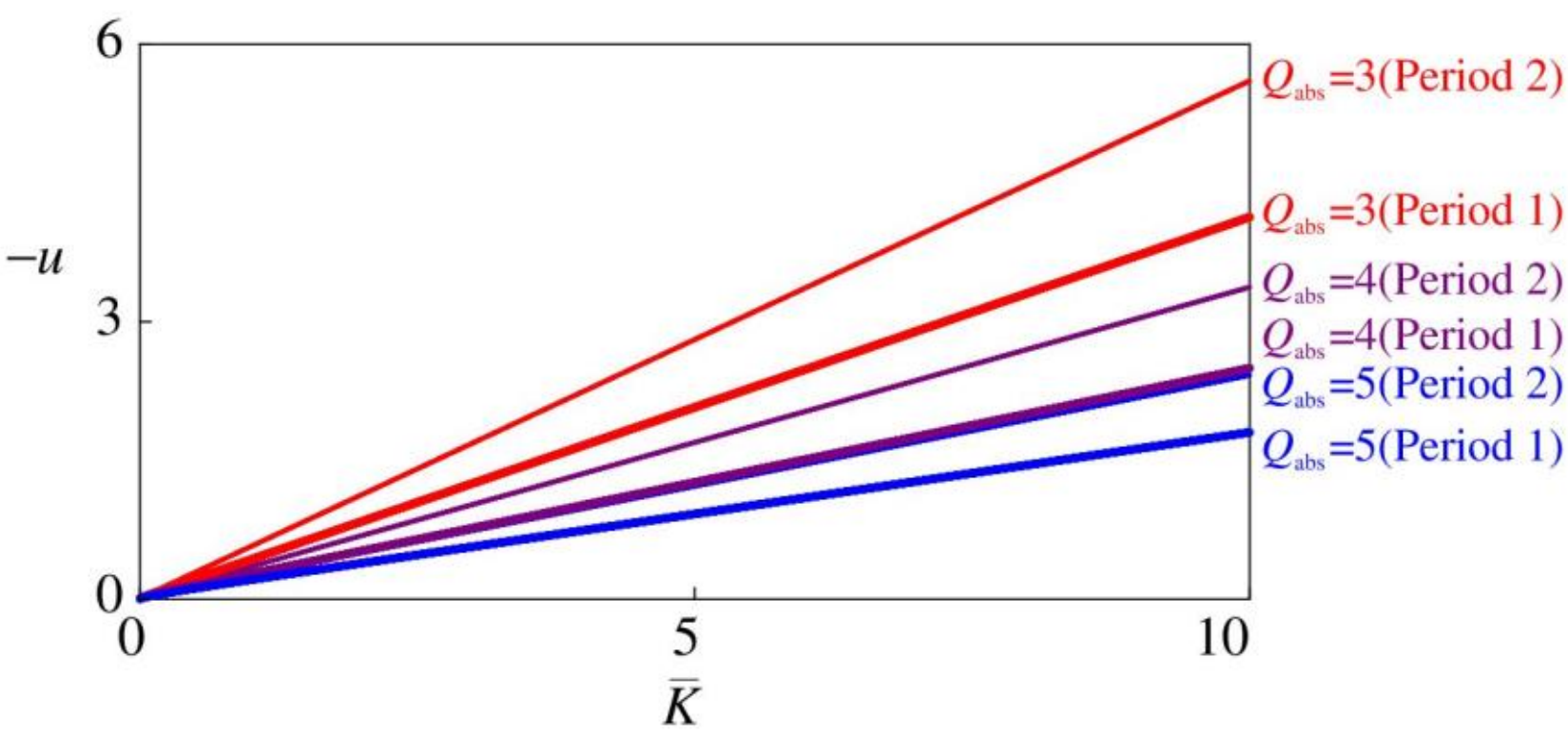

**Figure 9.** The computed optimal $-u$ with respect to the cost constraint $\overline{K}$.

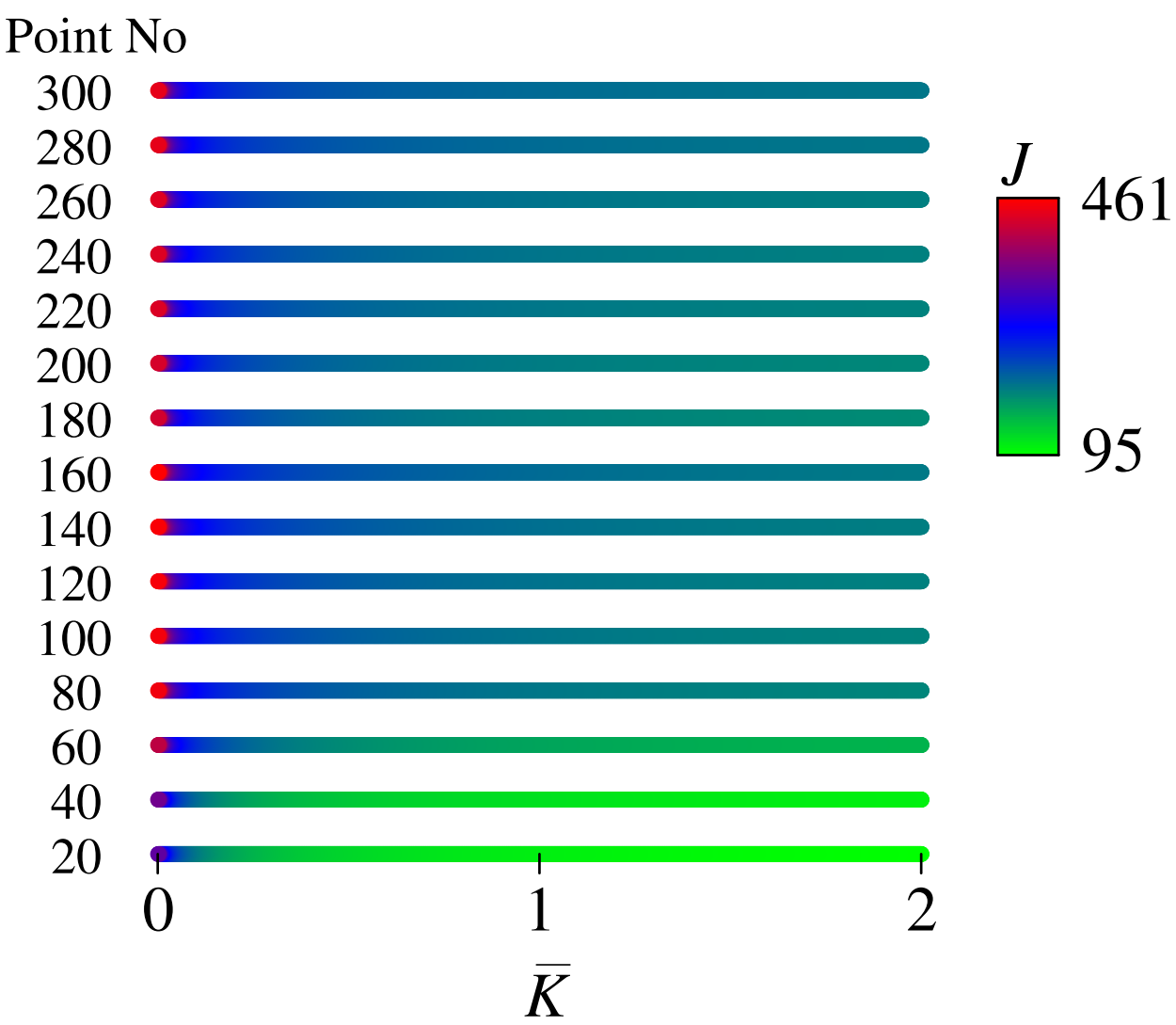

**Figure 10.** The computed minimum variance $J$ with respect to the cost constraint $\bar{K}$ for all the specific grid points. The target abstraction is fixed to $Q_{\mathrm{abs}} = 5$ ($\mathrm{m}^3$/s) in Period 1.

**Table 1.** Minimizing and maximizing locations of the minimized variance $J$ among all specific grid points. We examined the cost constraints $\bar{K}$ (0.001, 0.002, 0.003,…, 10.000 ($\mathrm{m}^6/\mathrm{s}^2$)) and the target abstractions $Q_{\mathrm{abs}}$ (1–8 ($\mathrm{m}^3$/s)). Minimizing point No. in Periods 1 and 2 are 20 in all the cases. Maximizing point No. in Period 2 is 300 in all the cases.

| $Q_{\mathrm{abs}}$ ($\mathrm{m}^3$/s) | Maximizing point No. in Period 1 |
|---|---|
| 1 | 160 |
| 2 | 160 |
| 3 | 160 $\left(\bar{K} \le 5.318\right)$ |
| | 300 (Otherwise) |
| 4 | 160 $\left(\bar{K} \le 0.977\right)$ |
| | 300 (Otherwise) |
| 5 | 160 $\left(\bar{K} \le 0.952\right)$ |
| | 300 (Otherwise) |
| 6 | 160 $\left(\bar{K} \le 1.152\right)$ |
| | 300 (Otherwise) |
| 7 | 160 $\left(\bar{K} \le 1.536\right)$ |
| | 300 (Otherwise) |
| 8 | 160 $\left(\bar{K} \le 2.200\right)$ |
| | 300 (Otherwise) |

### 4.5 Extended model

A cost-efficient control of the streamflow is vital, whereas strongly regulating the streamflow close to a fixed target value reduces the variability [49]. Herein, we considered an extended model with an additional constraint to maintain the streamflow variability. Indeed, flow regulations have been concerned as a potential stressor against local river environments and ecosystems. The problem is:

$$\text{Find } \inf_{\rho,u} J \text{ , where } J = \limsup_{T\to+\infty} \frac{1}{T}\mathbb{E}\left[\int_0^T \left(Q_t - \hat{Q}\right)^2 \mathrm{d}t\right] \tag{63}$$

subject to the cost constraints

$$K = \limsup_{T\to+\infty} \frac{1}{T}\mathbb{E}\left[\int_0^T c_t^2 \mathrm{d}t\right] \le \bar{K} \quad \text{and} \quad P = \limsup_{T\to+\infty} \frac{1}{T}\mathbb{E}\left[\int_0^T C_t^2 \mathrm{d}t\right] \le \bar{P} \tag{64}$$

and the same parameter constraints. The second inequality in (64) represents the constraint that maintains the flow variability measured in terms of the difference between the inflow and outflow (i.e., $C = X - Y$) based on a prescribed upper bound $\bar{P} > 0$. Namely, smaller $P$ suggests smaller modification of the river discharge.

Through analytical calculations, we obtain

$$P = \left(q-1\right)^2 \left\{ \left(\mathbb{E}\left[Y_{n,t}\right]\right)^2 + \frac{1}{2}\frac{AM_2}{D^2}\sum_{i=1}^{n}\frac{c_i}{r_i}\frac{h}{r_i D + h} \right\} \tag{65}$$

for the finite-dimensional case and

$$P = \left(q-1\right)^2 \left\{ \left(\mathbb{E}\left[Y_t\right]\right)^2 + \frac{1}{2}\frac{AM_2}{D^2}\int_0^{+\infty}\frac{1}{r}\frac{h}{rD + h}\pi\left(\mathrm{d}r\right) \right\} \tag{66}$$

for the infinite-dimensional case, where $h = \rho - u > 0$. These formulae can be effectively incorporated into the optimization problem by considering $P$ as an increasing function of $h$. Namely, the optimal $h$ is the largest $h > 0$ satisfying the two constraints (64).

Because of the functional shape of (66), we have the lower and upper bounds of $P = P(h)$ as

$$\left(q-1\right)^2\left(\mathbb{E}\left[Y_t\right]\right)^2 \le P\left(h\right) \le \left(q-1\right)^2 \mathbb{E}\left[\left(Y_t\right)^2\right], \tag{67}$$

showing that the problem is non-feasible if $\bar{P} < \left(q-1\right)^2\left(\mathbb{E}\left[Y_t\right]\right)^2$ and is trivial if $\left(q-1\right)^2 \mathbb{E}\left[\left(Y_t\right)^2\right]$. Thus, we assume

$$\left(q-1\right)^2\left(\mathbb{E}\left[Y_t\right]\right)^2 \le \bar{P} \le \left(q-1\right)^2 \mathbb{E}\left[\left(Y_t\right)^2\right]. \tag{68}$$

For demonstration, we set $\bar{P} = 100$ ($m^6/s^2$) in Period 1 as we have an upper bound of 167.8 ($m^6/s^2$) and a lower bound of 30.4 ($m^6/s^2$) based on (67). Similar results were obtained at the other points. The problem is, therefore, nontrivial, and the second constraint in (64) will be activated. **Figure 11** shows the computed minimized $J$ at Point 180 with the parameter region where the first and second constraints (64) were activated. The first constraint of (64), the cost constraint, is activated for small $\bar{K}$; otherwise, the second constraint on the flow variability is activated to maintain the flow variability within the required range. For relatively small water abstraction and cost constraints, and when the problem is subject to small $\bar{K}$, the additional constraint is not activated.

As demonstrated in this section, the proposed mathematical framework can be applied to the problem subject to multiple constraints without resorting to complex numerical algorithms owing to its high tractability.

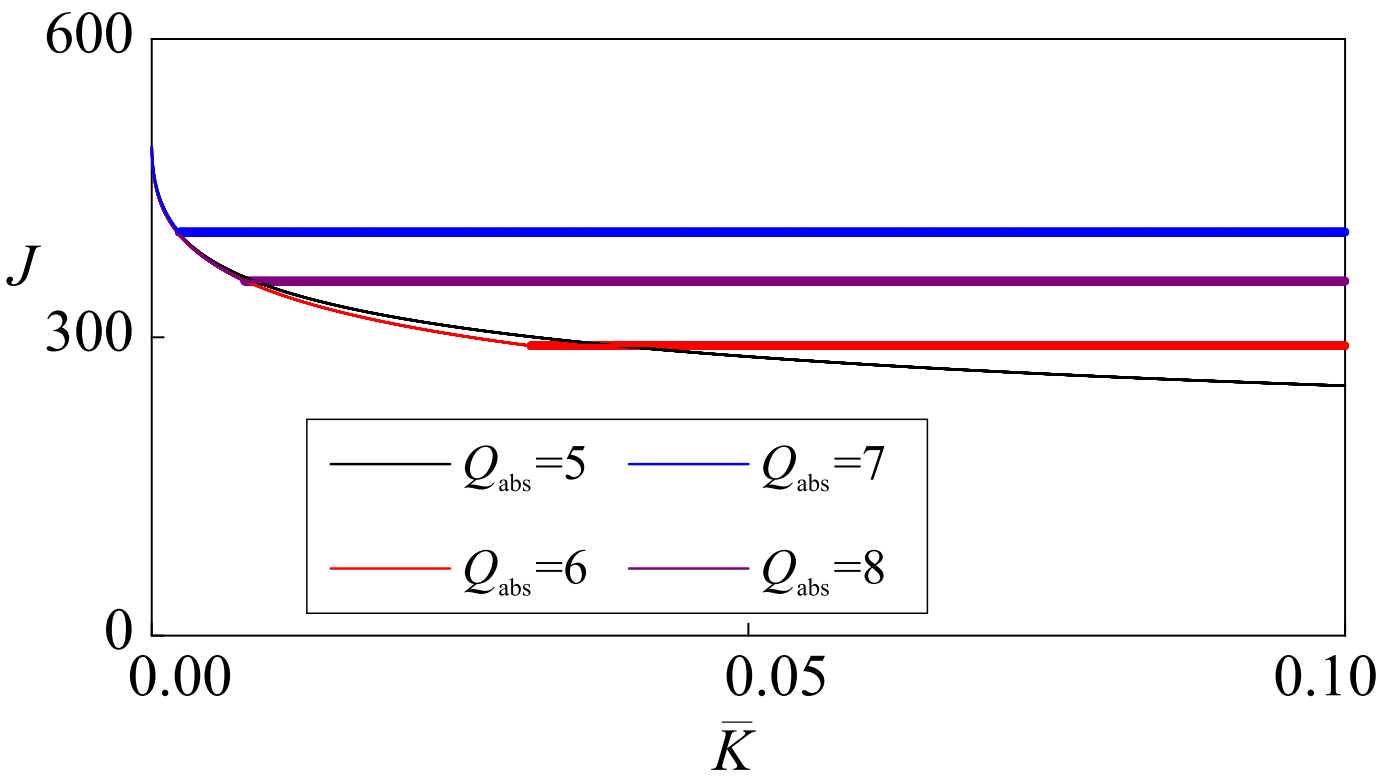

**Figure 11.** The computed minimized $J$ at Point 180 in Period 1 with the parameter region where the first and second constraints (64) are activated indicated by bold lines. The bold-line region is not appearing in the plot area for the case $Q_{abs} = 5$ ($m^3$/s).

## 5. Conclusion

We investigated a new optimization problem of river discharge based on the supCBI process to obtain the polynomially (subexponentially) decaying ACF. We found both theoretically and numerically that the infinite-dimensional nature of the optimization problem can be resolved by the Markovian lift. In addition, it was found that the optimization problem can be solved uniquely under certain conditions, and that an implementable numerical algorithm for solving the problem can be designed using a common methodology. The Markovian lift has been suggested to be a powerful mathematical tool for modelling and analysis of MMA processes. Our framework is therefore innovative in theory and not computationally demanding. The parameter values were identified from the data, and the model was analyzed.

The finding of this study recommends that modelling the streamflow dynamics should be based on a mathematical model that is able to capture the subexponential ACF. The supCBI process was successful in this case as demonstrated in the previous section. In addition, designing a water abstraction point based on the proposed model should consider the statistics of the discharge in target river environments because the objective and cost of the proposed model significantly depend on them. Water abstraction necessarily alters the streamflow dynamics, while they can be better managed by minimizing an objective like that proposed in this study.

The proposed model does not consider other environmental variables, such as the water temperature, suspended solid concentration, and some population density of aquatic fauna, which is a limitation. In addition, the seasonality of environmental and anthropogenic factors was not considered. In the studied river, local field managers especially concern water temperature and seasonality affecting growth and population dynamics of aquatic species. The affine properties of the proposed framework can be effectively employed in such applications if the system dynamics remain affine. In such cases, the corresponding static control problem would be computed numerically. Nevertheless, the computational complexity of BKE would not be prohibitive by invoking the affine properties, which allows the reduction of the equation to a system of ordinary differential equations. The model reduction can also be employed to reduce the problem size. Model identification is another potential issue to be tucked for approaching to more complex problems. A modelling and identification study concerning fish migration is currently undergoing.

For future studies, applications of the MMA process to more realistic stochastic control for water resource management would be investigated where the non-Markovian nature of the MMA process would be analyzed using a regression Monte Carlo method. In this case, the storage volume will be added to the stochastic system, whose value should be constrained in a compact set representing the dam capacity, resulting in a more complex and interesting problem. The model may contain uncertain parameter values possibly due to climate uncertainty, resulting in a worst-case optimization problem. Coupling with the price dynamics would also be investigated.

**Acknowledgments:** We used ACCMS, Kyoto University, Japan for a part of the computation.

**Funding:** This study was supported by the Japan Society for the Promotion of Science (22K14441, 22H02456), Environmental Research Projects from the Sumitomo Foundation (203160), and a Grant from MLIT Japan (B4R202002).

**Competing interests:** The author declares that there are no competing interests in this paper.

**References**

[1] L. Lu, Q. Tang, H. Li, Z. Li, Damming river shapes distinct patterns and processes of planktonic bacterial and microeukaryotic communities. Environ. Microbiol. 24 (2022) 1760-1774. https://doi.org/10.1111/1462-2920.15872.
[2] J.C. Stanton, J. Marek, L.S. Hall, B.E. Kus, A. Alvarado, B.K. Orr, E. Morrissette, L. Riege, W.E. Thogmartin, Recovery planning in a dynamic system: integrating uncertainty into a decision support tool for an endangered songbird. Ecol. Soc. 24 (2019) 11. https://doi.org/10.5751/ES-11169-240411.
[3] E. Akomeah, L.A. Morales-Marın, M. Carr, A. Sadeghian, K.E. Lindenschmidt, The impacts of changing climate and streamflow on nutrient speciation in a large Prairie reservoir. J. Environ. Manage. 288 (2021) 112262. https://doi.org/10.1016/j.jenvman.2021.112262.
[4] M. Carlin, M. Redolfi, M. Tubino, The long-term response of alternate bars to the hydrological regime. Water. Resour. Res. 57 (2021) e2020. https://doi.org/10.1029/2020WR029314.
[5] M.I. Brunner, L. Slater, L.M. Tallaksen, M. Clark, Challenges in modeling and predicting floods and droughts: a review. Wiley Interdiscip. Rev. Water. 8 (2021) e1520. https://doi.org/10.1002/wat2.1520.
[6] V.T.G. Boulomytis, A.C. Zuffo, M.A. Imteaz, Detection of flood influence criteria in ungauged basins on a combined Delphi-AHP approach. Oper. Res. Perspect. 6 (2019) 100116. https://doi.org/10.1016/j.orp.2019.100116.
[7] G. Bouveret, R. Dumitrescu, P. Tankov, Technological change in water use: A mean-field game approach to optimal investment timing. Oper. Res. Perspect. 9 (2022) 100225. https://doi.org/10.1016/j.orp.2022.100225.
[8] H. Li, P. Liu, S. Guo, Q. Zuo, J. Tao, K. Huang, Z. Yang, D. Han, B. Ming, Integrating teleconnection factors into long-term complementary operating rules for hybrid power systems: A case study of Longyangxia hydro-photovoltaic plant in China. Renew. Energy. 186 (2022) 517-534. https://doi.org/10.1016/j.renene.2022.01.034.
[9] G. Botter, S. Basso, I. Rodriguez-Iturbe, A. Rinaldo A. Resilience of river flow regimes. Proc. Natl. Acad. Sci. USA 110 (2013) 12925-12930. https://doi.org/10.1073/pnas.1311920110.
[10] M.B. Bertagni, P. Perona, C. Camporeale, Parametric transitions between bare and vegetated states in water-driven patterns. Proc. Natl. Acad. Sci, USA 115 (2018) 8125-8130. https://doi.org/10.1073/pnas.1721765115.
[11] A.J. Parolari, S. Pelrine, M.S. Bartlett, Stochastic water balance dynamics of passive and controlled stormwater basins. Adv. Water. Resour. 122 (2018) 328-339. https://doi.org/10.1016/j.advwatres.2018.10.016.
[12] A.B.V. Mello, M.C.S. Lima, A.D.C. Nascimento, A notable Gamma-Lindley first-order autoregressive process: an application to hydrological data. Environmetrics. 33 (2022) e2724. https://doi.org/10.1002/env.2724.
[13] N. Durighetto, V. Mariotto, F. Zanetti, K.J. McGuire, G. Mendicino, A Senatore, G. Botter, Probabilistic description of streamflow and active length regimes in rivers. Water. Resour. Res. (2022) e2021. https://doi.org/10.1029/2021WR031344.
[14] J.M. Turowski, Upscaling sediment-flux-dependent fluvial bedrock incision to long timescales. JGR Earth. Surface. (2021) e2020. https://doi.org/10.1029/2020JF005880.
[15] H. Yoshioka, Y. Yoshioka, Designing cost-efficient inspection schemes for stochastic streamflow environment using an effective Hamiltonian approach. Optim. Eng. 23 (2022) 1375-1407. https://doi.org/10.1007/s11081-021-09655-7.
[16] M.B. Al Sawaf, K. Kawanisi, Assessment of mountain river streamflow patterns and flood events using information and complexity measures. J. Hydrol. 590 (2020) 125508. https://doi.org/10.1016/j.jhydrol.2020.125508.
[17] H. Yoshioka, Fitting a superposition of Ornstein-Uhlenbeck processes to time series of discharge in a perennial river environment. ANZIAMJ. 63 (2022) C84-C96. https://doi.org/10.21914/anziamj.v63.16985.

[18] M. Girons Lopez, L. Crochemore, I.G. Pechlivanidis, Benchmarking an operational hydrological model for providing seasonal forecasts in Sweden. Hydrol. Earth. Syst. Sci. 25 (2021) 1189-1209. https://doi.org/10.5194/hess-25-1189-2021.
[19] C. Chen, Y. Tian, Y.K. Zhang, X. He, X. Yang, X. Liang, Y. Zheng, F. Han, C. Zheng, C. Yang, Effects of agricultural activities on the temporal variations of streamflow: trends and long memory. Stoch. Environ. Res. Risk. Assess. 33 (2019) 1553-1564. https://doi.org/10.1007/s00477-019-01714-x.
[20] L. Spezia, A. Vinten, R. Paroli, M. Stutter, An evolutionary Monte Carlo method for the analysis of turbidity high-frequency time series through Markov switching autoregressive models. Environmetrics. 32 (2021) e2695. https://doi.org/10.1002/env.2695.
[21] A. Habib, Exploring the physical interpretation of long-term memory in hydrology. Stoch. Environ. Res. Risk. Assess. 34 (2020) 2083-2091. https://doi.org/10.1007/s00477-020-01883-0.
[22] S.J. Sutanto, H.A.J. Van Lanen, Catchment memory explains hydrological drought forecast performance. Sci. Rep. 12 (2022) 2689. https://doi.org/10.1038/s41598-022-06553-5.
[23] L. Wang, W. Xia, Power-type derivatives for rough volatility with jumps. J. Futures. Markets. 42 (2022) 1369-1406. https://doi.org/10.1002/fut.22337.
[24] M. Moser, R. Stelzer, Functional regular variation of Lévy-driven multivariate mixed moving average processes. Extremes. 16 (2013) 351-382. https://doi.org/10.1007/s10687-012-0165-y.
[25] V. Courgeau, A.E.D. Veraart, Likelihood theory for the graph Ornstein-Uhlenbeck process. Stat. Inference. Stoch. Processes. 25 (2022) 227-60. https://doi.org/10.1007/s11203-021-09257-1.
[26] H. Yoshioka, M. Tsujimura, Derivation and computation of integro-Riccati equation for ergodic control of infinite-dimensional SDE, in: D. Groen, C. de Mulatier, M. Paszynski, V.V. Krzhizhanovskaya, J.J. Dongarra, P.M.A. Sloot (Eds.), Computational Science. Lect. Notes Comput. Sci. ICCS 2022. Vol. 13353. ICCS, Springer, Cham, 2022, pp. 577-588.
[27] C. Cuchiero, J. Teichmann, Markovian lifts of positive semidefinite affine Volterra-type processes. Decis, Econ, Fin. 42 (2019) 407-448. https://doi.org/10.1007/s10203-019-00268-5.
[28] E. Abi Jaber, E. Miller, H. Pham, Linear-quadratic control for a class of stochastic Volterra equations: solvability and approximation. Ann. Appl. Probab. 31 (2021) 2244-2274. https://doi.org/10.1214/20-AAP1645.
[29] I. Fatkhullin, B. Polyak, Optimizing static linear feedback: gradient method. SIAM J. Control. Optim. 59 (2021) 3887-3911. https://doi.org/10.1137/20M1329858.
[30] V. Blondel, J.N. Tsitsiklis, NP-hardness of some linear control design problems. SIAM J. Control. Optim. 35 (1997) 2118-2127. https://doi.org/10.1137/S0363012994272630.
[31] Y. Tang, Y. Zheng, N. Li, Analysis of the optimization landscape of linear quadratic gaussian (lqg) control. Learn. Dyn. Contr. 144 (2021) 599-610. http://proceedings.mlr.press/v144/tang21a.html.
[32] H. Yoshioka, A supCBI process with application to streamflow discharge and a model reduction, preprint (2022). https://arxiv.org/abs/2206.05923.
[33] B. Øksendal, A. Sulem, Applied Stochastic Control of Jump Diffusions, third ed., Springer, Cham 2019.
[34] H. Yoshioka, M. Tsujimura, T. Tanaka, Y. Yoshioka, A. Hashiguchi, Modeling and computation of an integral operator Riccati equation for an infinite-dimensional stochastic differential equation governing streamflow discharge, Comput. Math. Appl. 126 (2022) 115-148. https://doi.org/10.1016/j.camwa.2022.09.009.
[35] S. Göttlich, R. Korn, K. Lux, Optimal control of electricity input given an uncertain demand. Math. Methods. Oper. Res. 90 (2019) 301-328. https://doi.org/10.1007/s00186-019-00678-6.
[36] V. Fasen, C. Klüppelberg, Extremes of supOU processes, in: F.E. Benth, G. Nunno, T. Lindstrøm, B. Øksendal, T. Zhang (Eds.), Springer, Berlin, Heidelberg, 2007. pp. 339-359.
[37] P. Hulsman, T.A. Bogaard, H.H.G. Savenije, Rainfall-runoff modelling using river-stage time series in the absence of reliable discharge information: a case study in the semi-arid mara River basin. Hydrol. Earth. Syst. Sci. 22 (2018) 5081-5095. https://doi.org/10.5194/hess-22-5081-2018.
[38] H. McMillan, Linking hydrologic signatures to hydrologic processes: a review. Hydrol. Processes. 34 (2020) 1393-1409. https://doi.org/10.1002/hyp.13632.
[39] W.M. Wonham, On a matrix Riccati equation of stochastic control. SIAM J. Control. 6 (1968) 681-697. https://doi.org/10.1137/0306044.
[40] G. Choné, P.K. Biron, Assessing the relationship between river mobility and habitat. River. Res. Applic. 32 (2016) 528-539. https://doi.org/10.1002/rra.2896.
[41] A. Święch, J. Zabczyk, Integro-PDE in Hilbert Spaces: existence of viscosity solutions. Potential. Anal. 45 (2016) 703-736. https://doi.org/10.1007/s11118-016-9563-0.

[42] A. Klenke, Probability Theory: A Comprehensive Course, second ed., Springer, Cham, 2020.
[43] P. Greengard, V. Rokhlin, An algorithm for the evaluation of the incomplete gamma function. Adv. Comp. Math. 45 (2019) 23-49. https://doi.org/10.1007/s10444-018-9604-x.
[44] H. Yoshioka, T. Tanaka, F. Aranishi, M. Tsujimura, Y. Yoshioka, Impulsive fishery resource transporting strategies based on an open-ended stochastic growth model having a latent variable. Math. Methods. Appl. Sci. (2021). online published. https://doi.org/10.1002/mma.7982.
[45] S. Kitada, Long-term translocation explains population genetic structure of a recreationally fished iconic species in Japan: combining current knowledge with reanalysis. Aquaculture. Fish. Fish. 2 (2022) 130-145. https://doi.org/10.1002/aff2.34.
[46] A. Haddadchi, A. Kuczynski, J.T. Hoyle, C. Kilroy, D.J. Booker, M. Hicks, Periphyton removal flows determined by sediment entrainment thresholds. Ecol. Modell. 434 (2020) 109263. https://doi.org/10.1016/j.ecolmodel.2020.109263.
[47] H. Yoshioka, Y. Yoshioka, Stochastic streamflow and dissolved silica dynamics with application to the worst-case long-run evaluation of water environment. Optim. Eng. (2022). online published. https://doi.org/10.1007/s11081-022-09743-2.
[48] T. Tanaka, H. Yoshioka, Y. Yoshioka, DEM-based river cross-section extraction and 1-D streamflow simulation for eco-hydrological modeling: a case study in upstream Hiikawa River, Japan. Hydrol. Res. Lett. 15 (2021) 71-76. https://doi.org/10.3178/hrl.15.71.
[49] A.F. Van Loo, S. Rangecroft, G. Coxon, M. Werner, N. Wanders, G. Di Baldassarre, E. Tijdeman, M. Bosman, T. Gleeson, A. Nauditt, A. Aghakouchak, J.A. Breña-Naranjo, O. Cenobio-Cruz, A.C. Costa1, M. Fendekova, G. Jewitt, D.G. Kingston, J. Loft, S.M. Mager, I. Mallakpour, I. Masih, H. Maureira-Cortés, E. Toth, P. Van Oel, F. Van Ogtrop, K. Verbist, J.P. Vidal, L. Wen, M. Yu, X. Yuan, M. Zhang, Henny A.J. Van Lanen, Streamflow droughts aggravated by human activities despite management. Environ. Res. Lett. 17 (2022) 044059. https://doi.org/10.1088/1748-9326/ac5def.
[50] H. Yoshioka, M. Tsujimura, Hamilton–Jacobi–Bellman–Isaacs equation for rational inattention in the long-run management of river environments under uncertainty. Comput. Math. Appl. 112 (2022) 23-54. https://doi.org/10.1016/j.camwa.2022.02.013.

**Supplementary material of "Stochastic optimization of a mixed moving average process for controlling non-Markovian streamflow environments" by Hidekazu Yoshioka et al.**

## Appendix A

This appendix proves **Propositions 1-2**. The calculation is straightforward, while they contain several technical parts where we need to symmetrize some summands.

### A.1 Proof of Proposition 1

We guess $\Phi$ of the form

$$\Phi\left(x, y_1, ..., y_n\right) = \Phi\left(y_0, y_1, ..., y_n\right) = \sum_{i,j=0}^{n} \frac{1}{2} a_{i,j} y_i y_j + \sum_{i=0}^{n} b_i y_i \text{ , where } x = y_0 \tag{A.69}$$

with a symmetric matrix $a = \left[a_{i,j}\right]_{0 \le i,j \le n}$ and a vector $b = \left[b_i\right]_{0 \le i \le n}$. Then, we obtain

$$\frac{\partial \Phi}{\partial y_k}\left(x, y_1, ..., y_n\right) = \frac{1}{2}\sum_{i=0}^{n}\left(a_{i,k} + a_{k,i}\right) y_i + b_k = \sum_{i=0}^{n} a_{i,k} y_i + b_k \tag{A.70}$$

and

$$\begin{aligned}
&\int_0^{+\infty}\left(\Phi\left(x+z, ..., y_i + z, ...\right) - \Phi\right) v\left(\mathrm{d}z\right) \\
&= \int_0^{+\infty}\left(\begin{aligned}&\sum_{k,j=0}^{n}\frac{1}{2}a_{k,j}\left(y_k + \delta_{k,0}z + \delta_{k,i}z\right)\left(y_j + \delta_{j,0}z + \delta_{j,i}z\right) + \sum_{i=0}^{n} b_k\left(y_k + \delta_{k,0}z + \delta_{k,i}z\right) \\ &-\sum_{k,j=0}^{n}\frac{1}{2}a_{k,j}y_k y_j - \sum_{i=0}^{n} b_k y_k\end{aligned}\right) v\left(\mathrm{d}z\right) \\
&= \int_0^{+\infty}\left(\begin{aligned}&\sum_{k,j=0}^{n}\frac{1}{2}a_{k,j}\left(\delta_{k,0}\delta_{j,0} + \delta_{k,0}\delta_{j,i} + \delta_{k,i}\delta_{j,0} + \delta_{k,i}\delta_{j,i}\right)z^2 \\ &\sum_{k,j=0}^{n}\frac{1}{2}a_{k,j}\delta_{j,0}zy_k + \sum_{k,j=0}^{n}\frac{1}{2}a_{k,j}\delta_{j,i}zy_k + \sum_{k,j=0}^{n}\frac{1}{2}a_{k,j}\delta_{k,0}zy_j + \sum_{k,j=0}^{n}\frac{1}{2}a_{k,j}\delta_{k,i}zy_j \\ &+\left(b_0 + b_i\right)z\end{aligned}\right) v\left(\mathrm{d}z\right) \\
&= \int_0^{+\infty}\left(\frac{1}{2}\left(a_{0,0} + a_{0,i} + a_{i,0} + a_{i,i}\right)z^2 + \sum_{k=0}^{n}\frac{1}{2}\left(a_{k,0} + a_{k,i}\right)y_k z + \sum_{j=0}^{n}\frac{1}{2}\left(a_{0,j} + a_{i,j}\right)y_j z + \left(b_0 + b_i\right)z\right) v\left(\mathrm{d}z\right) \\
&= \frac{1}{2}M_2\left(a_{0,0} + a_{0,i} + a_{i,0} + a_{i,i}\right) + \frac{1}{2}M_1\sum_{k=0}^{n}\left(a_{0,k} + a_{k,0} + a_{k,i} + a_{i,k}\right)y_k + M_1\left(b_0 + b_i\right)
\end{aligned} \tag{A.71}$$

and hence

$$
\begin{aligned}
&\sum_{i=1}^{n}\left(c_i A+r_i B y_i\right)\int_0^{+\infty}\left(\Phi\left(x+z,...,y_i+z,...\right)-\Phi\right)v\left(\mathrm{d}z\right)\\
&=\sum_{i=1}^{n}\left(c_i A+r_i B y_i\right)\left\{\begin{array}{l}\frac{1}{2}M_2\left(a_{0,0}+a_{0,i}+a_{i,0}+a_{i,i}\right)\\ +\frac{1}{2}M_1\sum_{k=0}^{n}\left(a_{0,k}+a_{k,0}+a_{k,i}+a_{i,k}\right)y_k+M_1\left(b_0+b_i\right)\end{array}\right\}\\
&=\frac{1}{2}AM_2\sum_{i=1}^{n}c_i\left(a_{0,0}+a_{0,i}+a_{i,0}+a_{i,i}\right)+\frac{1}{2}BM_2\sum_{i=1}^{n}r_i\left(a_{0,0}+a_{0,i}+a_{i,0}+a_{i,i}\right)y_i\\
&+\frac{1}{2}AM_1\sum_{i=1}^{n}\sum_{k=0}^{n}c_i\left(a_{0,k}+a_{k,0}+a_{k,i}+a_{i,k}\right)y_k+\frac{1}{2}BM_1\sum_{i=1}^{n}\sum_{k=0}^{n}r_i\left(a_{0,k}+a_{k,0}+a_{k,i}+a_{i,k}\right)y_i y_k\\
&+AM_1\sum_{i=1}^{n}c_i\left(b_0+b_i\right)+BM_1\sum_{i=1}^{n}r_i y_i\left(b_0+b_i\right)
\end{aligned}
\quad, \quad (\mathrm{A}.72)
$$

$$
\begin{aligned}
\left\{-\left(\rho-u\right)x+\sum_{i=1}^{n}\left(\rho-r_i\right)y_i\right\}\frac{\partial\Phi}{\partial x}&=\left\{-\left(\rho-u\right)y_0+\sum_{i=1}^{n}\left(\rho-r_i\right)y_i\right\}\frac{\partial\Phi}{\partial y_0}\\
&=\left\{-\left(\rho-u\right)y_0+\sum_{i=1}^{n}\left(\rho-r_i\right)y_i\right\}\left(\sum_{k=0}^{n}\frac{1}{2}\left(a_{0,k}+a_{k,0}\right)y_k+b_0\right)\\
&=-\frac{1}{2}\left(\rho-u\right)\sum_{k=0}^{n}\left(a_{0,k}+a_{k,0}\right)y_0 y_k-\left(\rho-u\right)b_0 y_0\\
&\quad+\sum_{i=1}^{n}\sum_{k=0}^{n}\frac{1}{2}\left(\rho-r_i\right)\left(a_{0,k}+a_{k,0}\right)y_k y_i+b_0\sum_{i=1}^{n}\left(\rho-r_i\right)y_i
\end{aligned}
\quad, \quad (\mathrm{A}.73)
$$

and

$$
-\sum_{i=1}^{n}r_i y_i\frac{\partial\Phi}{\partial y_i}=-\frac{1}{2}\sum_{i=1}^{n}\sum_{k=0}^{n}r_i\left(a_{i,k}+a_{k,i}\right)y_i y_k-\sum_{i=1}^{n}r_i b_i y_i\ . \quad (\mathrm{A}.74)
$$

Consequently, the BKE is reduced to

$$
\begin{aligned}
&-J+\left(y_0\right)^2-2\hat{X}y_0+\hat{X}^2\\
&-\frac{1}{2}\left(\rho-u\right)\sum_{k=0}^{n}\left(a_{0,k}+a_{k,0}\right)y_0 y_k-\left(\rho-u\right)b_0 y_0+\sum_{i=1}^{n}\sum_{k=0}^{n}\frac{1}{2}\left(\rho-r_i\right)\left(a_{0,k}+a_{k,0}\right)y_k y_i+b_0\sum_{i=1}^{n}\left(\rho-r_i\right)y_i\\
&-\frac{1}{2}\sum_{i=1}^{n}\sum_{k=0}^{n}r_i\left(a_{i,k}+a_{k,i}\right)y_i y_k-\sum_{i=1}^{n}r_i b_i y_i\\
&+\frac{1}{2}AM_2\sum_{i=1}^{n}c_i\left(a_{0,0}+a_{0,i}+a_{i,0}+a_{i,i}\right)+\frac{1}{2}BM_2\sum_{i=1}^{n}r_i\left(a_{0,0}+a_{0,i}+a_{i,0}+a_{i,i}\right)y_i\\
&+\frac{1}{2}AM_1\sum_{i=1}^{n}\sum_{k=0}^{n}c_i\left(a_{0,k}+a_{k,0}+a_{k,i}+a_{i,k}\right)y_k+\frac{1}{2}BM_1\sum_{i=1}^{n}\sum_{k=0}^{n}r_i\left(a_{0,k}+a_{k,0}+a_{k,i}+a_{i,k}\right)y_i y_k\\
&+AM_1\sum_{i=1}^{n}c_i\left(b_0+b_i\right)+BM_1\sum_{i=1}^{n}r_i y_i\left(b_0+b_i\right)=0
\end{aligned}
$$

(A.75)

Using the symmetry of $a$, we can rewrite (A.75) as

$$
\begin{aligned}
&-J+\left(y_0\right)^2-2\hat{X}y_0+\hat{X}^2\\
&-\frac{1}{2}(\rho-u)\sum_{k=0}^{n}\left(a_{0,k}+a_{k,0}\right)y_0y_k-\left(\rho-u\right)b_0y_0\\
&+\frac{1}{4}\sum_{i=1}^{n}\sum_{k=0}^{n}\left\{\left(\rho-r_i\right)\left(a_{0,k}+a_{k,0}\right)+\left(\rho-r_k\right)\left(a_{0,i}+a_{i,0}\right)\right\}y_ky_i+b_0\sum_{i=1}^{n}\left(\rho-r_i\right)y_i\\
&-\frac{1}{4}\sum_{i=1}^{n}\sum_{k=0}^{n}\left\{r_i\left(a_{i,k}+a_{k,i}\right)+r_k\left(a_{k,i}+a_{i,k}\right)\right\}y_iy_k-\sum_{i=1}^{n}r_ib_iy_i\\
&+\frac{1}{2}AM_2\sum_{i=1}^{n}c_i\left(a_{0,0}+a_{0,i}+a_{i,0}+a_{i,i}\right)+\frac{1}{2}BM_2\sum_{i=1}^{n}r_i\left(a_{0,0}+a_{0,i}+a_{i,0}+a_{i,i}\right)y_i\\
&+\frac{1}{2}AM_1\sum_{i=1}^{n}\sum_{k=0}^{n}c_i\left(a_{0,k}+a_{k,0}+a_{k,i}+a_{i,k}\right)y_k+\frac{1}{4}BM_1\sum_{i=1}^{n}\sum_{k=0}^{n}\left\{\begin{matrix}r_i\left(a_{0,k}+a_{k,0}+a_{k,i}+a_{i,k}\right)\\+r_k\left(a_{0,i}+a_{i,0}+a_{i,k}+a_{k,i}\right)\end{matrix}\right\}y_iy_k\\
&+AM_1\sum_{i=1}^{n}c_i\left(b_0+b_i\right)+BM_1\sum_{i=1}^{n}r_iy_i\left(b_0+b_i\right)=0
\end{aligned}
\quad . \tag{A.76}
$$

Comparing each term of (A.76) leads to the following identities.

**(Coefficient proportional to 1)**

$$
-J+\hat{X}^2+\frac{1}{2}AM_2\sum_{i=1}^{n}c_i\left(a_{0,0}+2a_{0,i}+a_{i,i}\right)+AM_1\sum_{i=1}^{n}c_i\left(b_0+b_i\right)=0\,. \tag{A.77}
$$

**(Coefficient proportional to $y_0$)**

$$
-2\hat{X}-\left(\rho-u\right)b_0+AM_1\sum_{i=1}^{n}c_i\left(a_{0,0}+a_{0,i}\right)=0\,. \tag{A.78}
$$

**(Coefficient proportional to $y_i$, $1\le i\le n$)**

$$
b_0\left(\rho-r_i\right)-r_ib_i+\frac{1}{2}BM_2r_i\left(a_{0,0}+2a_{0,i}+a_{i,i}\right)+AM_1\sum_{k=1}^{n}c_k\left(a_{0,i}+a_{k,i}\right)+r_iBM_1\left(b_0+b_i\right)=0\,. \tag{A.79}
$$

**(Coefficient proportional to $\left(y_0\right)^2$)**

$$
1-\left(\rho-u\right)a_{0,0}=0\,. \tag{A.80}
$$

**(Coefficient proportional to $y_0y_i$, $1\le i\le n$)**

$$
-\left(\rho-u\right)a_{0,i}+\left(\rho-r_i\right)a_{0,0}-r_ia_{0,i}+BM_1r_i\left(a_{0,0}+a_{0,i}\right)=0\,. \tag{A.81}
$$

**(Coefficient proportional to $y_iy_j$, $1\le i,j\le n$)**

$$
\left(\rho-r_i\right)a_{0,j}-r_ia_{i,j}+BM_1r_i\left(a_{0,j}+a_{i,j}\right)+\left(\rho-r_j\right)a_{0,i}-r_ja_{i,j}+BM_1r_j\left(a_{0,i}+a_{i,j}\right)=0\,. \tag{A.82}
$$

Put $D=1-BM_1$, $p_i=\dfrac{r_i}{\rho-u}$, and $q=\dfrac{\rho}{\rho-u}$. From (A.80), we obtain

$$
a_{0,0}=\frac{1}{\rho-u}>0\,. \tag{A.83}
$$

Substituting (A.83) into (A.81) yields

$$a_{0,i} = \frac{1}{\rho - u} \frac{q - p_i D}{p_i D + 1}, \quad 1 \le i \le n. \tag{A.84}$$

By (A.82) and (A.84), we obtain

$$a_{i,j} = \frac{1}{\rho - u} \frac{(q - p_i D)(q - p_j D)}{(p_i + p_j) D} \left\{ \frac{1}{p_j D + 1} + \frac{1}{p_i D + 1} \right\}, \quad 1 \le i, j \le n. \tag{A.85}$$

Now, we have

$$a_{0,0} + a_{0,i} = \frac{1}{\rho - u} \frac{q + 1}{p_i D + 1}, \quad 1 \le i \le n. \tag{A.86}$$

Then, from (A.78) and (A.86), we obtain

$$b_0 = -\frac{2}{\rho - u} \hat{X} + \frac{1}{\rho - u} (q + 1) A M_1 \sum_{i=1}^{n} \frac{c_i}{r_i} \frac{p_i}{p_i D + 1}. \tag{A.87}$$

We have

$$\begin{aligned} a_{0,0} + 2a_{0,i} + a_{i,i} &= (a_{0,0} + a_{0,i}) + (a_{0,i} + a_{i,i}) \\ &= \frac{1}{\rho - u} \frac{q+1}{p_i D + 1} + \frac{1}{\rho - u} \frac{q - p_i D}{p_i D + 1} + \frac{1}{\rho - u} \frac{(q - p_i D)^2}{p_i D} \frac{1}{p_i D + 1}, \quad 1 \le i \le n. \\ &= \frac{1}{r_i} \frac{p_i D + q^2}{D(p_i D + 1)} \end{aligned} \tag{A.88}$$

From (A.79), (A.84), (A.85), (A.88), we obtain

$$\begin{aligned} & r_i (1 - B M_1)(b_i + b_0) \\ &= \rho b_0 + \frac{1}{2} B M_2 r_i (a_{0,0} + 2a_{0,i} + a_{i,i}) + A M_1 \sum_{k=1}^{n} c_k (a_{0,i} + a_{k,i}) \\ &= -2q\hat{X} + q(q+1) A M_1 \sum_{k=1}^{n} \frac{c_k}{r_k} \frac{p_k}{p_k D + 1} + \frac{1}{2} B M_2 \frac{p_i D + q^2}{D(p_i D + 1)} \qquad , \quad 1 \le i \le n. \\ &\quad + A M_1 \frac{q - p_i D}{p_i D + 1} \sum_{k=1}^{n} \frac{c_k}{r_k} p_k + A M_1 \sum_{k=1}^{n} \frac{c_k p_k}{r_k} \frac{(q - p_i D)(q - p_k D)}{(p_i + p_k) D} \left\{ \frac{1}{p_i D + 1} + \frac{1}{p_k D + 1} \right\} \end{aligned} \tag{A.89}$$

We can proceed as

$$\begin{aligned} & r_i D (b_i + b_0) \\ &= -2q\hat{X} + q(q+1) A M_1 \sum_{k=1}^{n} \frac{c_k}{r_k} \frac{p_k}{p_k D + 1} + \frac{1}{2} B M_2 \frac{p_i D + q^2}{D(p_i D + 1)} \qquad , \quad 1 \le i \le n. \\ &\quad + A M_1 \frac{q - p_i D}{p_i D + 1} \sum_{k=1}^{n} \frac{c_k}{r_k} p_k + A M_1 \sum_{k=1}^{n} \frac{c_k p_k}{r_k} \frac{(q - p_i D)(q - p_k D)}{(p_i + p_k) D} \left\{ \frac{1}{p_i D + 1} + \frac{1}{p_k D + 1} \right\} \end{aligned} \tag{A.90}$$

Then, we obtain

$$
\begin{aligned}
J &= \hat{X}^{2} + AM_1 \sum_{i=1}^{n} c_i \left( b_0 + b_i \right) + \frac{1}{2} AM_2 \frac{1}{\rho - u} \sum_{i=1}^{n} c_i \frac{p_i + q^2}{p_i \left( p_i + 1 \right)} \\
&= \hat{X}^{2} + AM_1 \sum_{i=1}^{n} \frac{c_i}{r_i D} \begin{pmatrix} -2q\hat{X} + q\left(q+1\right) AM_1 \sum_{k=1}^{n} \frac{c_k}{r_k} \frac{p_k}{p_k D + 1} + \frac{1}{2} BM_2 \frac{p_i D + q^2}{D\left(p_i D + 1\right)} \\ +AM_1 \frac{q - p_i D}{p_i D + 1} \sum_{k=1}^{n} \frac{c_k}{r_k} p_k \\ +AM_1 \sum_{k=1}^{n} \frac{c_k p_k}{r_k} \frac{\left(q - p_i D\right)\left(q - p_k D\right)}{\left(p_i + p_k\right) D} \left\{ \frac{1}{p_i D + 1} + \frac{1}{p_k D + 1} \right\} \end{pmatrix} \\
&\quad + \frac{1}{2} AM_2 \sum_{i=1}^{n} \frac{c_i}{r_i} \frac{p_i D + q^2}{D\left(p_i D + 1\right)} \\
&= \left( \hat{X} - \frac{qAM_1}{D} \sum_{i=1}^{n} \frac{c_i}{r_i} \right)^2 - \frac{q^2}{D^2} \left( AM_1 \right)^2 \sum_{i=1}^{n} \sum_{k=1}^{n} \frac{c_i}{r_i} \frac{c_k}{r_k} \\
&\quad + \frac{\left(AM_1\right)^2}{D} \sum_{i=1}^{n} \sum_{k=1}^{n} \frac{c_i}{r_i} \frac{c_k}{r_k} p_k \left\{ \begin{array}{l} \frac{q\left(q+1\right)}{p_k D + 1} + \frac{q - p_i D}{p_i D + 1} \\ + \frac{\left(q - p_i D\right)\left(q - p_k D\right)}{\left(p_i + p_k\right) D} \left\{ \frac{1}{p_i D + 1} + \frac{1}{p_k D + 1} \right\} \end{array} \right\} \\
&\quad + \frac{1}{2} \frac{AM_2}{D^2} \sum_{i=1}^{n} \frac{c_i}{r_i} \frac{p_i D + q^2}{p_i D + 1}
\end{aligned} \tag{A.91}
$$

We have

$$
\begin{aligned}
& p_k \left\{ \frac{q\left(q+1\right)}{p_k D + 1} + \frac{q - p_i D}{p_i D + 1} + \frac{\left(q - p_i D\right)\left(q - p_k D\right)}{\left(p_i + p_k\right) D} \left\{ \frac{1}{p_i D + 1} + \frac{1}{p_k D + 1} \right\} \right\} - \frac{q^2}{D} \\
&= p_k \frac{\left(q - p_i D\right)\left(p_i D + q\right)}{D\left(p_i D + 1\right)\left(p_i + p_k\right)} - p_i \frac{\left(q - p_k D\right)\left(p_k D + q\right)}{D\left(p_k D + 1\right)\left(p_i + p_k\right)}
\end{aligned} \tag{A.92}
$$

Due to the equation via a symmetry of (A.92), we obtain

$$
\begin{aligned}
& \sum_{i=1}^{n} \sum_{k=1}^{n} \frac{c_i}{r_i} \frac{c_k}{r_k} \left[ p_k \frac{\left(q - p_i D\right)\left(p_i D + q\right)}{D\left(p_i D + 1\right)\left(p_i + p_k\right)} - p_i \frac{\left(q - p_k D\right)\left(p_k D + q\right)}{D\left(p_k D + 1\right)\left(p_i + p_k\right)} \right] \\
&= \sum_{i=1}^{n} \sum_{k=1}^{n} \frac{c_i}{r_i} \frac{c_k}{r_k} \frac{1}{p_i + p_k} p_k \frac{\left(q - p_i D\right)\left(p_i D + q\right)}{D\left(p_i D + 1\right)} - \sum_{i=1}^{n} \sum_{k=1}^{n} \frac{c_i}{r_i} \frac{c_k}{r_k} \frac{1}{p_i + p_k} p_i \frac{\left(q - p_k D\right)\left(p_k D + q\right)}{D\left(p_k D + 1\right)} \\
&= \sum_{i=1}^{n} \sum_{k=1}^{n} \frac{c_i}{r_i} \frac{c_k}{r_k} \frac{1}{p_i + p_k} p_k \frac{\left(q - p_i D\right)\left(p_i D + q\right)}{D\left(p_i D + 1\right)} - \sum_{i=1}^{n} \sum_{k=1}^{n} \frac{c_k}{r_k} \frac{c_i}{r_i} \frac{1}{p_k + p_i} p_k \frac{\left(q - p_i D\right)\left(p_i D + q\right)}{D\left(p_i D + 1\right)} \\
&= \sum_{i=1}^{n} \sum_{k=1}^{n} \frac{c_i}{r_i} \frac{c_k}{r_k} \frac{1}{p_i + p_k} p_k \frac{\left(q - p_i D\right)\left(p_i D + q\right)}{D\left(p_i D + 1\right)} - \sum_{i=1}^{n} \sum_{k=1}^{n} \frac{c_i}{r_i} \frac{c_k}{r_k} \frac{1}{p_i + p_k} p_k \frac{\left(q - p_i D\right)\left(p_i D + q\right)}{D\left(p_i D + 1\right)} \\
&= 0
\end{aligned} \tag{A.93}
$$

Hence, we arrive at the desired result

$$
\begin{aligned}
J &= \left( \hat{X} - \frac{qAM_1}{D} \sum_{i=1}^{n} \frac{c_i}{r_i} \right)^2 + \frac{1}{2} \frac{AM_2}{D} \left( 1 + \frac{BM_1}{D} \right) \sum_{i=1}^{n} \frac{c_i}{r_i} \frac{p_i D + q^2}{p_i D + 1} \\
&= \left( \hat{X} - \frac{qAM_1}{D} \sum_{i=1}^{n} \frac{c_i}{r_i} \right)^2 + \frac{1}{2} \frac{AM_2}{D^2} \sum_{i=1}^{n} \frac{c_i}{r_i} \frac{p_i D + q^2}{p_i D + 1}
\end{aligned} \tag{A.94}
$$

combined with $r_i = p_i(\rho - u)$.

Finally, the fact that $J$ of (A.94) equals that in (15) follows from the verification argument similar to the argument in Appendix A.1 of Yoshioka and Tsujimura [50].

□

**A.2 Proof of Proposition 2**

The proof is similar to that in **Section A.1** but with slightly different coefficients due to dealing with the BKE having a different term from that for computing $J$. The proof is technical especially for some part dealing with multiple summations as detailed below.

We guess $\Psi$ of the form

$$\Psi(x, y_1, ..., y_n) = \Psi(y_0, y_1, ..., y_n) = \sum_{i,j=0}^{n} \frac{1}{2} a_{i,j} y_i y_j + \sum_{i=0}^{n} b_i y_i \text{ , where } x = y_0 . \tag{A.95}$$

with a symmetric matrix $a = [a_{i,j}]_{0 \le i,j \le n}$ and a vector $b = [b_i]_{0 \le i \le n}$. By straightforward calculations as in **Section A.1**, we obtain

$$\begin{aligned}
& -L + (y_0)^2 - 2q\sum_{i=1}^{n} y_0 y_i + q^2 \sum_{i=1}^{n}\sum_{k=1}^{n} y_i y_k \\
& -\frac{1}{2}(\rho - u)\sum_{k=0}^{n}(a_{0,k} + a_{k,0}) y_0 y_k - (\rho - u) b_0 y_0 \\
& +\sum_{i=1}^{n}\sum_{k=0}^{n}\frac{1}{2}(\rho - r_i)(a_{0,k} + a_{k,0}) y_k y_i + b_0 \sum_{i=1}^{n}(\rho - r_i) y_i \\
& -\frac{1}{2}\sum_{i=1}^{n}\sum_{k=0}^{n} r_i (a_{i,k} + a_{k,i}) y_i y_k - \sum_{i=1}^{n} r_i b_i y_i \\
& +\frac{1}{2} AM_2 \sum_{i=1}^{n} c_i (a_{0,0} + a_{0,i} + a_{i,0} + a_{i,i}) + \frac{1}{2} BM_2 \sum_{i=1}^{n} r_i (a_{0,0} + a_{0,i} + a_{i,0} + a_{i,i}) y_i \\
& +\frac{1}{2} AM_1 \sum_{i=1}^{n}\sum_{k=0}^{n} c_i (a_{0,k} + a_{k,0} + a_{k,i} + a_{i,k}) y_k + \frac{1}{2} BM_1 \sum_{i=1}^{n}\sum_{k=0}^{n} r_i (a_{0,k} + a_{k,0} + a_{k,i} + a_{i,k}) y_i y_k \\
& + AM_1 \sum_{i=1}^{n} c_i (b_0 + b_i) + BM_1 \sum_{i=1}^{n} r_i y_i (b_0 + b_i) = 0
\end{aligned} \quad . \tag{A.96}$$

With the symmetrization, we obtain the following identities.

**(Coefficient proportional to 1)**

$$-\bar{C} + \frac{1}{2} AM_2 \sum_{i=1}^{n} c_i (a_{0,0} + 2a_{0,i} + a_{i,i}) + AM_1 \sum_{i=1}^{n} c_i (b_0 + b_i) = 0 . \tag{A.97}$$

**(Coefficient proportional to $y_0$)**

$$-(\rho - u) b_0 + AM_1 \sum_{i=1}^{n} c_i (a_{0,0} + a_{0,i}) = 0 . \tag{A.98}$$

**(Coefficient proportional to $y_i$, $1 \le i \le n$)**

$$b_0\left(\rho-r_i\right)-r_ib_i+\frac{1}{2}BM_2r_i\left(a_{0,0}+2a_{0,i}+a_{i,i}\right)+AM_1\sum_{k=1}^{n}c_k\left(a_{0,i}+a_{k,i}\right)+r_iBM_1\left(b_0+b_i\right)=0\,. \quad \text{(A.99)}$$

**(Coefficient proportional to $\left(y_0\right)^2$)**

$$1-\left(\rho-u\right)a_{0,0}=0\,. \quad \text{(A.100)}$$

**(Coefficient proportional to $y_0y_i$, $1\le i\le n$)**

$$-2q-\left(\rho-u\right)a_{0,i}+\left(\rho-r_i\right)a_{0,0}-r_ia_{0,i}+BM_1r_i\left(a_{0,0}+a_{0,i}\right)=0\,. \quad \text{(A.101)}$$

**(Coefficient proportional to $y_iy_j$, $1\le i,j\le n$)**

$$2q^2+\left(\rho-r_i\right)a_{0,j}-r_ia_{i,j}+BM_1r_i\left(a_{0,j}+a_{i,j}\right)+\left(\rho-r_j\right)a_{0,i}-r_ja_{i,j}+BM_1r_j\left(a_{0,i}+a_{i,j}\right)=0\,. \quad \text{(A.102)}$$

By solving these equations, we obtain

$$a_{0,0}=\frac{1}{\rho-u}>0\,, \quad \text{(A.103)}$$

$$a_{0,i}=-\frac{1}{\rho-u}\frac{q+p_iD}{p_iD+1}\,,\ \ 1\le i\le n\,, \quad \text{(A.104)}$$

$$a_{i,j}=\frac{1}{\rho-u}\frac{1}{\left(p_i+p_j\right)D}\left\{-\frac{\left(q-p_iD\right)\left(q+p_jD\right)}{p_jD+1}-\frac{\left(q-p_jD\right)\left(q+p_iD\right)}{p_iD+1}+2q^2\right\}\,,\qquad 1\le i,j\le n\,, \quad \text{(A.105)}$$

$$a_{0,0}+a_{0,i}=-\frac{1}{\rho-u}\frac{q-1}{p_iD+1}\,,\ \ 1\le i\le n\,, \quad \text{(A.106)}$$

$$b_0=-AM_1\frac{q-1}{\rho-u}\sum_{i=1}^{n}\frac{c_i}{r_i}\frac{p_i}{p_iD+1}\,, \quad \text{(A.107)}$$

$$a_{0,0}+2a_{0,i}+a_{i,i}=\frac{1}{\rho-u}\frac{-\left(q+p_iD\right)^2+\left(p_iD+q^2\right)\left(p_iD+1\right)}{p_iD\left(p_iD+1\right)}\,,\ \ 1\le i\le n\,, \quad \text{(A.108)}$$

and hence

$$\begin{aligned}
&r_i\left(1-BM_1\right)\left(b_i+b_0\right)\\
&=-q\left(q-1\right)AM_1\sum_{k=1}^{n}\frac{c_k}{r_k}\frac{p_k}{p_kD+1}+\frac{1}{2}BM_2\frac{-\left(q+p_iD\right)^2+\left(p_iD+q^2\right)\left(p_iD+1\right)}{D\left(p_iD+1\right)}\\
&-AM_1\frac{q+p_iD}{p_iD+1}\sum_{k=1}^{n}\frac{c_k}{r_k}p_k\\
&+AM_1\sum_{k=1}^{n}\frac{c_k}{r_k}p_k\frac{1}{\left(p_i+p_k\right)D}\left\{-\frac{\left(q-p_iD\right)\left(q+p_kD\right)}{p_kD+1}-\frac{\left(q-p_kD\right)\left(q+p_iD\right)}{p_iD+1}+2q^2\right\}
\end{aligned}\,,\qquad 1\le i\le n\,. \quad \text{(A.109)}$$

We can proceed as

$$
\begin{aligned}
L &= \frac{1}{2} AM_2 \sum_{i=1}^{n} c_i \left( a_{0,0} + 2a_{0,i} + a_{i,i} \right) + AM_1 \sum_{i=1}^{n} c_i \left( b_0 + b_i \right) \\
&= \frac{1}{2} AM_2 \sum_{i=1}^{n} c_i \, \frac{1}{\rho - u} \frac{-\left(q + p_i D\right)^2 + \left(p_i D + q^2\right)\left(p_i D + 1\right)}{p_i D \left(p_i D + 1\right)} \\
&\quad + AM_1 \sum_{i=1}^{n} \frac{c_i}{r_i D} \left\{ \begin{array}{l}
-q\left(q-1\right) AM_1 \displaystyle\sum_{k=1}^{n} \frac{c_k}{r_k} \frac{p_k}{p_k D + 1} \\
+ \dfrac{1}{2} BM_2 \dfrac{-\left(q + p_i D\right)^2 + \left(p_i D + \dfrac{q^2}{2}\right)\left(p_i D + 1\right)}{D\left(p_i D + 1\right)} \\
- AM_1 \dfrac{q + p_i D}{p_i D + 1} \displaystyle\sum_{k=1}^{n} \frac{c_k}{r_k} p_k \\
+ AM_1 \displaystyle\sum_{k=1}^{n} \frac{c_k}{r_k} p_k \frac{1}{\left(p_i + p_k\right) D} \left\{ \begin{array}{l}
- \dfrac{\left(q - p_i D\right)\left(q + p_k D\right)}{p_k D + 1} \\
- \dfrac{\left(q - p_k D\right)\left(q + p_i D\right)}{p_i D + 1} + 2q^2
\end{array} \right\}
\end{array} \right\} \\
&= \frac{1}{2} \frac{AM_2}{D^2} \sum_{i=1}^{n} \frac{c_i}{r_i} \frac{-\left(q + p_i D\right)^2 + \left(p_i D + q^2\right)\left(p_i D + 1\right)}{p_i D + 1} \\
&\quad + \frac{\left(AM_1\right)^2}{D} \sum_{i=1}^{n} \sum_{k=1}^{n} \frac{c_i}{r_i} \frac{c_k}{r_k} p_k \left\{ \begin{array}{l}
- \dfrac{q\left(q-1\right)}{p_k D + 1} - \dfrac{q + p_i D}{p_i D + 1} \\
+ \dfrac{1}{\left(p_i + p_k\right) D} \left\{ \begin{array}{l}
- \dfrac{\left(q - p_i D\right)\left(q + p_k D\right)}{p_k D + 1} \\
- \dfrac{\left(q - p_k D\right)\left(q + p_i D\right)}{p_i D + 1} + 2q^2
\end{array} \right\}
\end{array} \right\} . \qquad \text{(A.110)}
\end{aligned}
$$

Now, we have

$$
\begin{aligned}
&\sum_{i=1}^{n}\sum_{k=1}^{n}\frac{c_i}{r_i}\frac{c_k}{r_k}p_k\left\{-\frac{q\left(q-1\right)}{p_kD+1}-\frac{q+p_iD}{p_iD+1}+\frac{1}{\left(p_i+p_k\right)D}\left\{\begin{array}{l}-\dfrac{\left(q-p_iD\right)\left(q+p_kD\right)}{p_kD+1}\\ -\dfrac{\left(q-p_kD\right)\left(q+p_iD\right)}{p_iD+1}+2q^2\end{array}\right\}\right\}\\
&=\sum_{i=1}^{n}\sum_{k=1}^{n}\frac{c_i}{r_i}\frac{c_k}{r_k}p_k\left\{\begin{array}{l}-\dfrac{q+p_iD}{p_iD+1}\\ -\dfrac{1}{\left(p_i+p_k\right)D}\dfrac{\left(q-p_iD\right)\left(q+p_kD\right)}{p_kD+1}\\ -\dfrac{1}{\left(p_i+p_k\right)D}\dfrac{\left(q-p_kD\right)\left(q+p_iD\right)}{p_iD+1}\end{array}\right\}\\
&-q\left(q-1\right)\sum_{i=1}^{n}\sum_{k=1}^{n}\frac{c_i}{r_i}\frac{c_k}{r_k}\frac{p_k}{p_kD+1}+2\sum_{i=1}^{n}\sum_{k=1}^{n}\frac{c_i}{r_i}\frac{c_k}{r_k}\frac{p_kq^2}{\left(p_i+p_k\right)D}\\
&=\sum_{i=1}^{n}\sum_{k=1}^{n}\frac{c_i}{r_i}\frac{c_k}{r_k}\frac{p_k}{\left(p_i+p_k\right)D}\left\{\begin{array}{l}-\dfrac{\left(q+p_iD\right)\left(p_i+p_k\right)D}{p_iD+1}\\ -\dfrac{\left(q-p_iD\right)\left(q+p_kD\right)}{p_kD+1}\\ -\dfrac{\left(q-p_kD\right)\left(q+p_iD\right)}{p_iD+1}\end{array}\right\}\\
&-\frac{q\left(q-1\right)}{2}\sum_{i=1}^{n}\sum_{k=1}^{n}\frac{c_i}{r_i}\frac{c_k}{r_k}\left(\frac{p_k}{p_kD+1}+\frac{p_i}{p_iD+1}\right)+\sum_{i=1}^{n}\sum_{k=1}^{n}\frac{c_i}{r_i}\frac{c_k}{r_k}\frac{\left(p_i+p_k\right)q^2}{\left(p_i+p_k\right)D}\\
&=-\sum_{i=1}^{n}\sum_{k=1}^{n}\frac{c_i}{r_i}\frac{c_k}{r_k}\frac{p_k}{\left(p_i+p_k\right)D}\left\{\frac{\left(q-p_iD\right)\left(q+p_kD\right)}{p_kD+1}+\frac{\left(q+p_iD\right)^2}{p_iD+1}\right\}\\
&-\frac{q\left(q-1\right)}{2}\sum_{i=1}^{n}\sum_{k=1}^{n}\frac{c_i}{r_i}\frac{c_k}{r_k}\left(\frac{p_k}{p_kD+1}+\frac{p_i}{p_iD+1}\right)+\frac{q^2}{D}\left(\sum_{i=1}^{n}\frac{c_i}{r_i}\right)^2 \quad . \qquad \text{(A.111)}
\end{aligned}
$$

We can proceed as

$$
\begin{aligned}
&\sum_{i=1}^{n}\sum_{k=1}^{n}\frac{c_i}{r_i}\frac{c_k}{r_k}\frac{p_k}{\left(p_i+p_k\right)D}\left\{\frac{\left(q-p_iD\right)\left(q+p_kD\right)}{p_kD+1}+\frac{\left(q+p_iD\right)^2}{p_iD+1}\right\} \\
&=\sum_{i=1}^{n}\sum_{k=1}^{n}\frac{c_i}{r_i}\frac{c_k}{r_k}\frac{p_k}{\left(p_i+p_k\right)D}\frac{\left(q-p_iD\right)\left(q+p_kD\right)}{p_kD+1}+\sum_{i=1}^{n}\sum_{k=1}^{n}\frac{c_i}{r_i}\frac{c_k}{r_k}\frac{p_k}{\left(p_i+p_k\right)D}\frac{\left(q+p_iD\right)^2}{p_iD+1} \\
&=\sum_{i=1}^{n}\sum_{k=1}^{n}\frac{c_i}{r_i}\frac{c_k}{r_k}\frac{p_k}{\left(p_i+p_k\right)D}\frac{\left(q-p_iD\right)\left(q+p_kD\right)}{p_kD+1}+\sum_{k=1}^{n}\sum_{i=1}^{n}\frac{c_k}{r_k}\frac{c_i}{r_i}\frac{p_i}{\left(p_i+p_k\right)D}\frac{\left(q+p_kD\right)^2}{p_kD+1} \\
&=\sum_{i=1}^{n}\sum_{k=1}^{n}\frac{c_i}{r_i}\frac{c_k}{r_k}\frac{1}{\left(p_i+p_k\right)D}\left\{\frac{p_k\left(q-p_iD\right)\left(q+p_kD\right)}{p_kD+1}+\frac{p_i\left(q+p_kD\right)^2}{p_kD+1}\right\} \\
&=\sum_{i=1}^{n}\sum_{k=1}^{n}\frac{c_i}{r_i}\frac{c_k}{r_k}\frac{q+p_kD}{\left(p_i+p_k\right)\left(p_kD+1\right)D}\left\{p_k\left(q-p_iD\right)+p_i\left(q+p_kD\right)\right\} \\
&=\sum_{i=1}^{n}\sum_{k=1}^{n}\frac{c_i}{r_i}\frac{c_k}{r_k}\frac{q+p_kD}{\left(p_i+p_k\right)\left(p_kD+1\right)D}q\left(p_i+p_k\right) \\
&=\sum_{i=1}^{n}\sum_{k=1}^{n}\frac{c_i}{r_i}\frac{c_k}{r_k}\frac{q\left(q+p_kD\right)}{\left(p_kD+1\right)D} \\
&=\frac{q}{2D}\sum_{i=1}^{n}\sum_{k=1}^{n}\frac{c_i}{r_i}\frac{c_k}{r_k}\left(\frac{q+p_iD}{p_iD+1}+\frac{q+p_kD}{p_kD+1}\right)
\end{aligned}
\quad . \qquad \text{(A.112)}
$$

Consequently, we obtain

$$
\begin{aligned}
L=&\frac{1}{2}\frac{AM_2}{D^2}\sum_{i=1}^{n}\frac{c_i}{r_i}\frac{-\left(q+p_iD\right)^2+\left(p_iD+q^2\right)\left(p_iD+1\right)}{p_iD+1} \\
&+\frac{\left(AM_1\right)^2q}{2D^2}\sum_{i=1}^{n}\sum_{k=1}^{n}\frac{c_i}{r_i}\frac{c_k}{r_k}\left\{-\left(\frac{q+p_iD}{p_iD+1}+\frac{q+p_kD}{p_kD+1}\right)-\left(q-1\right)D\left(\frac{p_k}{p_kD+1}+\frac{p_i}{p_iD+1}\right)+2q\right\}
\end{aligned}
\quad . \qquad \text{(A.113)}
$$

We have

$$
\sum_{i=1}^{n}\frac{c_i}{r_i}\frac{-\left(q+p_iD\right)^2+\left(p_iD+q^2\right)\left(p_iD+1\right)}{p_iD+1}=\sum_{i=1}^{n}\frac{c_i}{r_i}\frac{p_iD\left(-2q+1+q^2\right)}{p_iD+1}=\left(1-q\right)^2\frac{AM_2}{D^2}\sum_{i=1}^{n}\frac{c_i}{r_i}\frac{p_iD}{p_iD+1} \qquad \text{(A.114)}
$$

and

$$
\begin{aligned}
&-\left(\frac{q+p_iD}{p_iD+1}+\frac{q+p_kD}{p_kD+1}\right)-\left(q-1\right)D\left(\frac{p_k}{p_kD+1}+\frac{p_i}{p_iD+1}\right)+2q \\
&=-\left(\frac{q+p_iD+qDp_i-Dp_i}{p_iD+1}+\frac{q+p_kD+qDp_k+-p_k}{p_kD+1}\right)+2q \\
&=-2q+2q \\
&=0
\end{aligned}
\quad . \qquad \text{(A.115)}
$$

Substituting (A.114) and (A.115) into (A.113) yields the desired result

$$
L=\frac{1}{2}\left(1-q\right)^2\frac{AM_2}{D^2}\sum_{i=1}^{n}\frac{c_i}{r_i}\frac{p_iD}{p_iD+1}, \qquad \text{(A.116)}
$$

from which we obtain $K=\left(\rho-u\right)^2L$ by $r_i=p_i\left(\rho-u\right)$. Finally, the fact that $K$ presented here truly equals that in (16) follows from again by the verification argument similar to the argument in Appendix A.1 of Yoshioka and Tsujimura [50].

□

**A.3 Derivation of** $P$

By the definition of $P$, we have

$$P = \overline{\mathbb{E}\left[\left(Y_t - X_t\right)^2\right]} = \overline{\mathbb{E}\left[\left(X_t\right)^2\right]} - 2\overline{\mathbb{E}\left[X_t Y_t\right]} + \overline{\mathbb{E}\left[\left(Y_t\right)^2\right]}. \tag{A.117}$$

Here, for simplicity, we formally write $\overline{\mathbb{E}[\bullet]} = \lim_{T\to+\infty}\frac{1}{T}\mathbb{E}\left[\int_0^T(\bullet)\mathrm{d}t\right]$. We already have

$$J = \overline{\mathbb{E}\left[\left(X_t - \hat{X}\right)^2\right]} = \overline{\mathbb{E}\left[\left(X_t\right)^2\right]} - 2\hat{X}\,\overline{\mathbb{E}\left[X_t\right]} + \hat{X}^2 \tag{A.118}$$

or equivalently

$$\overline{\mathbb{E}\left[\left(X_t\right)^2\right]} = J + 2\hat{X}\,\overline{\mathbb{E}\left[X_t\right]} - \hat{X}^2 \tag{A.119}$$

and

$$\frac{C}{\left(\rho-u\right)^2} = \bar{C} = \overline{\mathbb{E}\left[\left(X_t - qY_t\right)^2\right]} = \overline{\mathbb{E}\left[\left(X_t\right)^2\right]} - 2q\overline{\mathbb{E}\left[X_t Y_t\right]} + q^2\,\overline{\mathbb{E}\left[\left(Y_t\right)^2\right]} \tag{A.120}$$

or equivalently by (A.119),

$$\overline{\mathbb{E}\left[X_t Y_t\right]} = \frac{1}{2q}\left\{\overline{\mathbb{E}\left[\left(X_t\right)^2\right]} - \bar{C} + q^2\,\overline{\mathbb{E}\left[\left(Y_t\right)^2\right]}\right\} = \frac{1}{2q}\left\{J - \hat{X}^2 - \bar{C} + 2\hat{X}\,\overline{\mathbb{E}\left[X_t\right]} + q^2\,\overline{\mathbb{E}\left[\left(Y_t\right)^2\right]}\right\}. \tag{A.121}$$

We then obtain

$$\begin{aligned}
P &= \overline{\mathbb{E}\left[\left(Y_t - X_t\right)^2\right]} \\
&= J + 2\hat{X}\,\overline{\mathbb{E}\left[X_t\right]} - \hat{X}^2 + \overline{\mathbb{E}\left[\left(Y_t\right)^2\right]} \\
&\quad -2\left\{\frac{1}{2q}\left\{J - \hat{X}^2 - \bar{C} + 2\hat{X}\,\overline{\mathbb{E}\left[X_t\right]} + q^2\,\overline{\mathbb{E}\left[\left(Y_t\right)^2\right]}\right\}\right\} \\
&= \left(1 - \frac{1}{q}\right)\left(J + 2\hat{X}\,\overline{\mathbb{E}\left[X_t\right]} - \hat{X}^2\right) + \frac{1}{q}\bar{C} + \left(1-q\right)\overline{\mathbb{E}\left[\left(Y_t\right)^2\right]}
\end{aligned}. \tag{A.122}$$

As $\overline{\mathbb{E}\left[X_t\right]} = \hat{X}$, we obtain

$$
\begin{aligned}
P &= \left(1-\frac{1}{q}\right)\left(J+\hat{X}^{2}\right)+\frac{1}{q}\bar{C}+(1-q)\overline{\mathbb{E}\left[\left(Y_t\right)^{2}\right]} \\
&= \frac{q-1}{q}\left(J+\hat{X}^{2}\right)+\frac{1}{q}\bar{C}+(1-q)\overline{\left\{\mathbb{V}\left[Y_t\right]+\left(\mathbb{E}\left[Y_t\right]\right)^{2}\right\}} \\
&= \frac{q-1}{q}\left(J+\hat{X}^{2}\right)+\frac{1}{q}\bar{C}+(1-q)\overline{\left\{\mathbb{V}\left[Y_t\right]+\left(\frac{\hat{X}}{q}\right)^{2}\right\}} \\
&= \frac{q-1}{q}\left(\frac{1}{2}\frac{AM_2}{D^2}\sum_{i=1}^{n}\frac{c_i}{r_i}\frac{p_iD+q^2}{p_iD+1}+\hat{X}^{2}\right)+\frac{1}{2q}(1-q)^2\frac{AM_2}{D^2}\sum_{i=1}^{n}\frac{c_i}{r_i}\frac{p_iD}{p_iD+1} \\
&\quad +(1-q)\left\{\frac{1}{2}\frac{AM_2}{D^2}\sum_{i=1}^{n}\frac{c_i}{r_i}+\left(\frac{\hat{X}}{q}\right)^{2}\right\} \\
&= \frac{q-1}{q}\hat{X}^{2}+\frac{1-q}{q^2}\hat{X}^{2}+\frac{1}{2}\frac{AM_2}{D^2}\left\{(1-q)\sum_{i=1}^{n}\frac{c_i}{r_i}+\frac{q-1}{q}\sum_{i=1}^{n}\frac{c_i}{r_i}\frac{p_iD+q^2}{p_iD+1}+\frac{1}{q}(1-q)^2\sum_{i=1}^{n}\frac{c_i}{r_i}\frac{p_iD}{p_iD+1}\right\} \\
&= \left(\frac{q-1}{q}\right)^{2}\hat{X}^{2}+\frac{1}{2}\frac{AM_2}{D^2}\frac{1-q}{q}\sum_{i=1}^{n}\frac{c_i}{r_i}\frac{q(1-q)}{p_iD+1}
\end{aligned}
$$

(A.123)

and consequently

$$
P=(q-1)^2\left\{\frac{\hat{X}^2}{q^2}+\frac{1}{2}\frac{AM_2}{D^2}\sum_{i=1}^{n}\frac{c_i}{r_i}\frac{1}{p_iD+1}\right\}. \tag{A.124}
$$

We obtain the desired result by combining (26) and (A.124). The proof for the infinite-dimensional case is essentially the same where each summation is replaced by a corresponding integration.

□

**Appendix B**

This appendix presents the convergence results of the integral $R=\int_0^{+\infty} r^{-1}\pi(\mathrm{d}r)$ with the Gamma density $\pi(\mathrm{d}r)=r^{\alpha-1}\exp(-r)\mathrm{d}r$, where the scaling parameter has been chosen to be 1 without any loss of generality. Here, we set $R_n=\sum_{i=1}^{n}c_i r_i^{-1}$. **Tables B.1-B.3** show that the convergence rate is around 0.5 for the shape parameters $\alpha = 1.8,\ 2.0,\ 2.2$. The convergence rate becomes lower as $\alpha$ decreases with which the regularity of the Gamma density near the origin becomes lower.

**Table B.1.** Convergence results of the integral $R=\int_0^{+\infty} r^{-1}\pi(\mathrm{d}r)$ using the Markovian lift: $\alpha=1.8$.

| $n$ | $R_n$ | $R$ | Relative error $e_n=(R-R_n)/R$ | Convergence rate $\log_2(e_{n-1}/e_n)$ |
|---|---|---|---|---|
| 64 | 1.15537 | 1.25 | 0.075704 | |
| 128 | 1.18043 | 1.25 | 0.055656 | 0.444 |
| 256 | 1.19886 | 1.25 | 0.040912 | 0.444 |
| 512 | 1.21242 | 1.25 | 0.030064 | 0.444 |
| 1024 | 1.22238 | 1.25 | 0.022096 | 0.444 |
| 2048 | 1.2297 | 1.25 | 0.016240 | 0.444 |
| 4096 | 1.23508 | 1.25 | 0.011936 | 0.444 |
| 8192 | 1.23904 | 1.25 | 0.008768 | 0.445 |

**Table B.2.** Convergence results of the integral $R=\int_0^{+\infty} r^{-1}\pi(\mathrm{d}r)$ using the Markovian lift: $\alpha=2.0$.

| $n$ | $R_n$ | $R$ | Relative error $e_n=(R-R_n)/R$ | Convergence rate $\log_2(e_{n-1}/e_n)$ |
|---|---|---|---|---|
| 64 | 0.94661 | 1 | 0.053390 | |
| 128 | 0.962226 | 1 | 0.037774 | 0.499 |
| 256 | 0.973281 | 1 | 0.026719 | 0.500 |
| 512 | 0.981103 | 1 | 0.018897 | 0.500 |
| 1024 | 0.986636 | 1 | 0.013364 | 0.500 |
| 2048 | 0.99055 | 1 | 0.009450 | 0.500 |
| 4096 | 0.993317 | 1 | 0.006683 | 0.500 |
| 8192 | 0.995274 | 1 | 0.004726 | 0.500 |

**Table B.3.** Convergence results of the integral $R=\int_0^{+\infty} r^{-1}\pi(\mathrm{d}r)$ using the Markovian lift: $\alpha=2.2$.

| $n$ | $R_n$ | $R$ | Relative error $e_n=(R-R_n)/R$ | Convergence rate $\log_2(e_{n-1}/e_n)$ |
|---|---|---|---|---|
| 64 | 0.800163 | 0.833333 | 0.039804 | |
| 128 | 0.810588 | 0.833333 | 0.027294 | 0.544 |
| 256 | 0.817741 | 0.833333 | 0.018710 | 0.545 |
| 512 | 0.822647 | 0.833333 | 0.012823 | 0.545 |
| 1024 | 0.82601 | 0.833333 | 0.008788 | 0.545 |
| 2048 | 0.828315 | 0.833333 | 0.006022 | 0.545 |
| 4096 | 0.829895 | 0.833333 | 0.004126 | 0.546 |
| 8192 | 0.830977 | 0.833333 | 0.002827 | 0.545 |

**Appendix C**

We present the model identification results mentioned in **Section 4.3** of the main text.

**Table C.1.** Identification results for Period 1 at three representative points.

| | Distance from O Dam | | | | | |
|---|---|---|---|---|---|---|
| Parameter | 2.3 (km): Point 20 | | 4.3 (km): Point 60 | | 10.2 (km): Point 180 | |
| $c_1$ (-) | 7.720.E−01 | | 8.113.E−01 | | 8.127.E−01 | |
| $c_2$ (s/m$^3$) | 4.434.E−03 | | 3.709.E−03 | | 3.942.E−03 | |
| $A$ (m$^{3c_1}$/s$^{c_1}$/h) | 2.391.E−02 | | 2.799.E−02 | | 2.918.E−02 | |
| $B$ (m$^{3(c_1-1)}$/s$^{c_1-1}$) | 3.637.E−02 | | 3.562.E−02 | | 3.603.E−02 | |
| $\underline{Q}$ (m$^3$/s) | 1.174.E+00 | | 1.226.E+00 | | 1.290.E+00 | |
| $\alpha$ (-) | 2.329.E+00 | | 2.248.E+00 | | 2.189.E+00 | |
| $D\beta$ (1/h) | 3.149.E−02 | | 3.328.E−02 | | 3.210.E−02 | |
| Statistics | Empirical | Model | Empirical | Model | Empirical | Model |
| Average (m$^3$/s) | 9.111.E+00 | 9.029.E+00 | 1.077.E+01 | 1.069.E+01 | 1.199.E+01 | 1.190.E+01 |
| Variance (m$^6$/s$^2$) | 4.023.E+02 | 4.039.E+02 | 4.796.E+02 | 4.813.E+02 | 5.022.E+02 | 5.040.E+02 |
| Skewness (-) | 1.305.E+01 | 1.175.E+01 | 1.339.E+01 | 1.206.E+01 | 1.228.E+01 | 1.106.E+01 |
| Kurtosis (-) | 2.430.E+02 | 2.526.E+02 | 2.619.E+02 | 2.723.E+02 | 2.205.E+02 | 2.293.E+02 |

**Table C.2.** Identification results for Period 2 at three representative points.

| | Distance from O Dam | | | | | |
|---|---|---|---|---|---|---|
| Parameter | 2.3 (km): Point 20 | | 4.3 (km): Point 60 | | 10.2 (km): Point 180 | |
| $c_1$ (-) | 8.400.E−01 | | 8.379.E−01 | | 8.381.E−01 | |
| $c_2$ (s/m$^3$) | 4.345.E−03 | | 3.645.E−03 | | 3.295.E−03 | |
| $A$ (m$^{3c_1}$/s$^{c_1}$/h) | 2.615.E−02 | | 2.907.E−02 | | 3.170.E−02 | |
| $B$ (m$^{3(c_1-1)}$/s$^{c_1-1}$) | 3.604.E−02 | | 3.511.E−02 | | 3.453.E−02 | |
| $\underline{Q}$ (m$^3$/s) | 1.836.E+00 | | 2.617.E+00 | | 2.752.E+00 | |
| $\alpha$ (-) | 1.865.E+00 | | 1.865.E+00 | | 1.874.E+00 | |
| $D\beta$ (1/h) | 2.941.E−02 | | 3.175.E−02 | | 3.226.E−02 | |
| Statistics | Empirical | Model | Empirical | Model | Empirical | Model |
| Average (m$^3$/s) | 1.622.E+01 | 1.610.E+01 | 1.784.E+01 | 1.769.E+01 | 1.920.E+01 | 1.903.E+01 |
| Variance (m$^6$/s$^2$) | 5.237.E+02 | 5.254.E+02 | 6.678.E+02 | 6.701.E+02 | 7.969.E+02 | 7.998.E+02 |
| Skewness (-) | 1.042.E+01 | 9.373.E+00 | 1.111.E+01 | 9.927.E+00 | 1.133.E+01 | 1.005.E+01 |
| Kurtosis (-) | 1.611.E+02 | 1.676.E+02 | 1.801.E+02 | 1.878.E+02 | 1.843.E+02 | 1.925.E+02 |

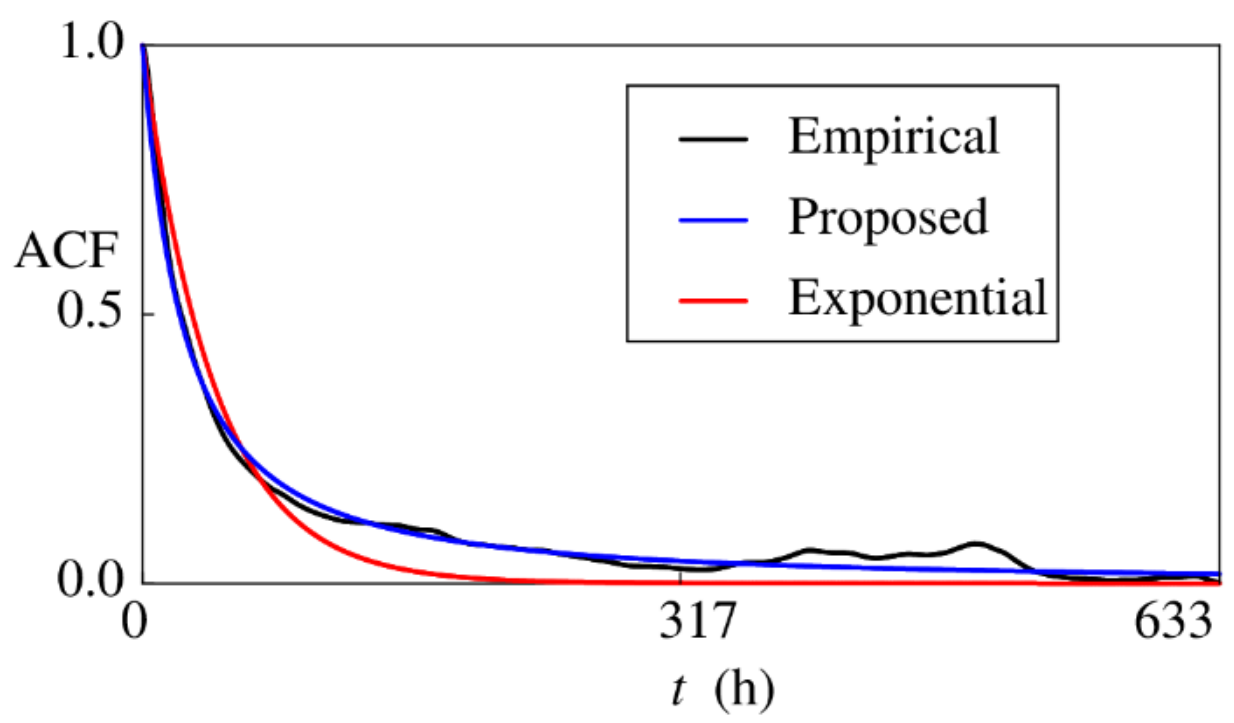

**Figure C.1.** Comparison of empirical and modelled ACFs of the discharge for Period 1 at Point 20.

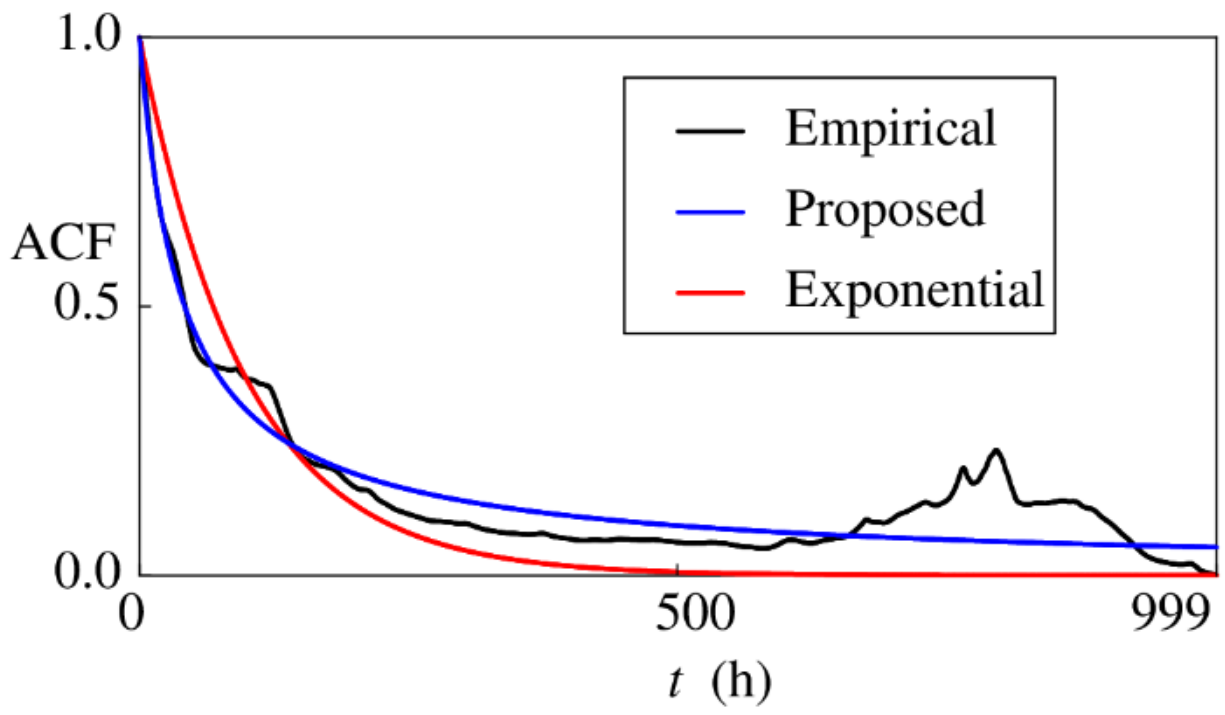

**Figure C.2.** Comparison of empirical and modelled ACFs of the discharge for Period 2 at Point 20.